\documentclass[10pt,a4paper]{article}
\usepackage{amsmath,amssymb,amsthm,graphics,graphicx,epic,eepic,multicol,ascmac}

\usepackage[all]{xy}

\numberwithin{equation}{section}
\theoremstyle{theorem}
\newtheorem{thm}{Theorem}[section]
\newtheorem{prop}[thm]{Proposition}
\newtheorem{lem}[thm]{Lemma}
\newtheorem{rem}[thm]{Remark}

\newtheorem{ex}[thm]{Example}

\theoremstyle{definition}
\newtheorem{defn}[thm]{Definition}

\def\al{\alpha}

\def\wht(#1){\widehat{\ #1\ }}

\newcommand{\cA}{{\mathcal A}}

\newcommand{\cF}{{\mathcal F}}

\newcommand{\frg}{\mathfrak g}
\newcommand{\frh}{\mathfrak h}

\newcommand{\frn}{\mathfrak n}
\newcommand{\frp}{\mathfrak p}

\newcommand{\frs}{\mathfrak s}

\newcommand{\bbC}{\mathbb C}

\newcommand{\ch}{\mathrm{ch}}

\newcommand{\lbr}{\begin{bmatrix}}
\newcommand{\rbr}{\end{bmatrix}}

\newcommand{\cd}{commutative diagram }

\def\ge{\frg}

\def\cP{{\mathcal P}}
\def\al{\alpha}

\def\beneme{\begin{enumerate}}
\def\beq{\begin{equation}}
\def\beqn{\begin{eqnarray}}
\def\beqnn{\begin{eqnarray*}}
\def\bfi{{\mathbf i}}
\def\bfii0{{\bf i_0}}

\def\bbra#1,#2,#3{\left\{\begin{array}{c}\hspace{-5pt}
#1;#2\\ \hspace{-5pt}#3\end{array}\hspace{-5pt}\right\}}
\def\cd{\cdots}
\def\ci(#1,#2){c_{#1}^{(#2)}}
\def\Ci(#1,#2){C_{#1}^{(#2)}}
\def\mpp(#1,#2,#3){#1^{(#2)}_{#3}}
\def\bCi(#1,#2){\ovl C_{#1}^{(#2)}}
\def\ch(#1,#2){c_{#2,#1}^{-h_{#1}}}
\def\cc(#1,#2){c_{#2,#1}}

\def\del{\delta}
\def\Del{\Delta}

\def\di(#1,#2){D_{#1}^{(#2)}}
\def\dbi(#1,#2){\ovl D_{#1}^{(#2)}}

\def\eneme{\end{enumerate}}

\def\eeq{\end{equation}}
\def\eeqn{\end{eqnarray}}
\def\eeqnn{\end{eqnarray*}}

\def\gau#1,#2{\left[\begin{array}{c}\hspace{-5pt}#1\\
\hspace{-5pt}#2\end{array}\hspace{-5pt}\right]}

\def\ji(#1,#2){j_{#1}^{(#2)}}

\def\lan{\langle}

\def\lm{\lambda}
\def\Lm{\Lambda}

\def\nd{\noindent}

\def\ovl{\overline}

\def\qq{\qquad}
\def\q{\quad}
\def\qed{\hfill\framebox[2mm]{}}

\def\ran{\rangle}

\def\TY(#1,#2,#3){#1^{(#2)}_{#3}}

\def\UU{{\mathcal U}}

\def\xxi(#1,#2,#3){\displaystyle {}^{#1}\Xi^{(#2)}_{#3}}
\def\xsi(#1,#2,#3){\displaystyle {}^{#1}\Sigma^{(#2)}_{#3}}
\def\xE(#1,#2,#3){\displaystyle {}^{#1}E_{#2}[#3]}
\def\xF(#1,#2){\displaystyle {}^{#1}F_{#2}}
\def\xx(#1,#2){\displaystyle {}^{#1}\Xi_{#2}}
\def\W1{W(\varpi_1)}

\def\m@th{\mathsurround=0pt}
\def\fsquare(#1,#2){
\hbox{\vrule$\hskip-0.4pt\vcenter to #1{\normalbaselines\m@th
\hrule\vfil\hbox to #1{\hfill$\scriptstyle #2$\hfill}\vfil\hrule}$\hskip-0.4pt
\vrule}}

\newcommand{\ba}{\begin{array}}
\newcommand{\ea}{\end{array}}

\newcommand{\eq}{\begin{eqnarray}}
\newcommand{\eneq}{\end{eqnarray}}

\def\m@th{\mathsurround=0pt}

\def\addsquare(#1,#2){\hbox{$
	\dimen1=#1 \advance\dimen1 by -0.8pt
	\vcenter to #1{\hrule height0.4pt depth0.0pt%
	\hbox to #1{%
        \vbox to \dimen1{\vss%
	\hbox to \dimen1{\hss$\scriptstyle~#2~$\hss}%
	\vss}%
	\vrule width0.4pt}%
	\hrule height0.4pt depth0.0pt}$}}

\def\Addsquare(#1,#2){\hbox{$
	\dimen1=#1 \advance\dimen1 by -0.8pt
	\vcenter to #1{\hrule height0.4pt depth0.0pt%
	\hbox to #1{%
	\vbox to \dimen1{\vss%
	\hbox to \dimen1{\hss$\scriptstyle~#2~$\hss}%
	\vss}%
	\vrule width0.4pt}%
	\hrule height0.4pt depth0.0pt}$}}

\def\fsquare(#1,#2){
\hbox{\vrule$\hskip-0.4pt\vcenter to #1{\normalbaselines\m@th
\hrule\vfil\hbox to #1{\hfill$\scriptstyle #2$\hfill}\vfil\hrule}$\hskip-0.4pt
\vrule}}

\def\seudosquare(#1,#2,#3){
\hbox{$\hskip-0.4pt\vcenter to #1{\normalbaselines\m@th
\hrule\vfil\hbox to #2{$\hfill\scriptstyle #3\hfill$}\vfil\hrule}$\hskip-0.4pt
\vrule}}

\def\topseudosquare(#1,#2,#3){
\hbox{$\hskip-0.4pt\vtop to #1{\normalbaselines\m@th
\hrule\vfil\hbox to #2{$\hfill\scriptstyle #3\hfill$}\vfil\hrule}$\hskip-0.4pt
\vrule}}

\def\Topseudosquare<#1,#2,#3>{
\hbox{$\hskip-0.4pt\vtop to #1{\normalbaselines\m@th
\hrule\vfil\hbox to #2{$\hfill\scriptstyle #3\hfill$}\vfil\hrule}$\hskip-0.4pt
\vrule}}

\def\TTopseudosquare<#1,#2,#3>{
\hbox{$\hskip-0.4pt\vtop to #1{\normalbaselines\m@th
\hrule\vfil\hbox to #2{$\hfill #3\hfill$}\vfil\hrule}$\hskip-0.4pt
\vrule}}

\def\tophalfsquare(#1,#2,#3){
\hbox{\vrule$\hskip-0.4pt\vtop to #1{\normalbaselines\m@th
\hrule\vfil\hbox to #2{\hfill$\scriptstyle #3$\hfill}\vfil\hrule}$\hskip-1pt
\vrule}}

\def\Tophalfsquare<#1,#2,#3>{
\hbox{\vrule$\hskip-0.4pt\vtop to #1{\normalbaselines\m@th
\hrule\vfil\hbox to #2{\hfill$\scriptstyle #3$\hfill}\vfil\hrule}$\hskip-1pt
\vrule}}

\def\TTophalfsquare<#1,#2,#3>{
\hbox{\vrule$\hskip-0.4pt\vtop to #1{\normalbaselines\m@th
\hrule\vfil\hbox to #2{\hfill$#3$\hfill}\vfil\hrule}$\hskip-1pt
\vrule}}

\def\threetwo(#1,#2,#3,#4,#5,#6){
\vcenter{\hbox{
        $\TTophalfsquare<1cm,1.0cm,#1>
         \TTopseudosquare<1cm,1.0cm,#2>$}
\vskip-1.2pt
\hbox{
        $\TTophalfsquare<1cm,1.0cm,#3>
         \TTopseudosquare<1cm,1.0cm,#4>$}
\vskip-1.2pt
\hbox{
        $\TTophalfsquare<1cm,1.0cm,#5>
         \TTopseudosquare<1cm,1.0cm,#6>$}
}}

\def\fourfourA(#1,#2,#3,#4,#5){
\vcenter{\hbox{
        $\tophalfsquare(0.8cm,1.2cm,#1)
         \topseudosquare(0.8cm,1.2cm,#1)
         \topseudosquare(0.8cm,1.2cm,#4)
         \topseudosquare(0.8cm,1.2cm,#1)$}
\vskip-1.2pt
\hbox{
        $\tophalfsquare(0.8cm,1.2cm,#2)
         \topseudosquare(0.8cm,1.2cm,#2)
         \topseudosquare(0.8cm,1.2cm,#4)
         \topseudosquare(0.8cm,1.2cm,#2)$}
\vskip-1.2pt
\hbox{
        $\tophalfsquare(1.0cm,1.2cm,#5)
         \topseudosquare(1.0cm,1.2cm,#5)
         \topseudosquare(1.0cm,1.2cm,#5)
         \topseudosquare(1.0cm,1.2cm,#5)$}
\vskip-1.4pt
\hbox{
        $\tophalfsquare(0.8cm,1.2cm,#3)
         \topseudosquare(0.8cm,1.2cm,#3)
         \topseudosquare(0.8cm,1.2cm,#4)
         \topseudosquare(0.8cm,1.2cm,#3)$}
}}

\def\fourfourC(#1,#2,#3,#4,#5){
\vcenter{\hbox{
        $\tophalfsquare(0.8cm,0.8cm,#1)
         \topseudosquare(0.8cm,0.8cm,#1)
         \topseudosquare(0.8cm,0.8cm,#4)
         \topseudosquare(0.8cm,0.8cm,#1)$}
\vskip-1.2pt
\hbox{
        $\tophalfsquare(0.8cm,0.8cm,#2)
         \topseudosquare(0.8cm,0.8cm,#2)
         \topseudosquare(0.8cm,0.8cm,#4)
         \topseudosquare(0.8cm,0.8cm,#2)$}
\vskip-1.2pt
\hbox{
        $\tophalfsquare(1.0cm,0.8cm,#5)
         \topseudosquare(1.0cm,0.8cm,#5)
         \topseudosquare(1.0cm,0.8cm,#5)
         \topseudosquare(1.0cm,0.8cm,#5)$}
\vskip-1.4pt
\hbox{
        $\tophalfsquare(0.8cm,0.8cm,#3)
         \topseudosquare(0.8cm,0.8cm,#3)
         \topseudosquare(0.8cm,0.8cm,#4)
         \topseudosquare(0.8cm,0.8cm,#3)$}
}}

\def\fourfourB(#1,#2,#3,#4,#5,#6,#7,#8){
\vcenter{\hbox{
        $\Tophalfsquare<1cm,2.0cm,#1>
         \topseudosquare(1cm,2.0cm,#7)
         \Topseudosquare<1cm,2.0cm,#2>$}
\vskip-1.5pt
\hbox{
        $\Tophalfsquare<1cm,2.0cm,#3>
         \topseudosquare(1cm,2.0cm,#7)
         \Topseudosquare<1cm,2.0cm,#4>$}
\vskip-1.5pt
\hbox{
        $\Tophalfsquare<1.0cm,2.0cm,#8>
         \topseudosquare(1.0cm,2.0cm,#8)
         \Topseudosquare<1.0cm,2.0cm,#8>$}
\vskip-1.4pt
\hbox{
        $\Tophalfsquare<1.0cm,2.0cm,#5>
         \topseudosquare(1.0cm,2.0cm,#7)
         \Topseudosquare<1.0cm,2.0cm,#6>$}
}}

\def\threefourA(#1,#2,#3,#4,#5,#6,#7,#8){
\vcenter{\hbox{
        $\Tophalfsquare<0.8cm,3.0cm,#1>
         \topseudosquare(0.8cm,3.0cm,#8)
         \topseudosquare(0.8cm,3.0cm,#2)
         \Topseudosquare<0.8cm,3.0cm,#3>$}
\vskip-1.5pt
\hbox{
        $\Tophalfsquare<0.8cm,3.0cm,#4>
         \topseudosquare(0.8cm,3.0cm,#8)
         \topseudosquare(0.8cm,3.0cm,#5)
         \Topseudosquare<0.8cm,3.0cm,#6>$}
\vskip-1.4pt
\hbox{
        $\Tophalfsquare<0.8cm,3.0cm,#7>
         \topseudosquare(0.8cm,3.0cm,#8)
         \topseudosquare(0.8cm,3.0cm,#7)
         \Topseudosquare<0.8cm,3.0cm,#7>$}
}}

\def\threefourB(#1,#2,#3,#4,#5,#6,#7,#8){
\vcenter{\hbox{
        $\Tophalfsquare<0.8cm,0.8cm,#1>
         \topseudosquare(0.8cm,0.8cm,#8)
         \topseudosquare(0.8cm,0.8cm,#2)
         \Topseudosquare<0.8cm,0.8cm,#3>$}
\vskip-1.5pt
\hbox{
        $\Tophalfsquare<0.8cm,0.8cm,#4>
         \topseudosquare(0.8cm,0.8cm,#8)
         \topseudosquare(0.8cm,0.8cm,#5)
         \Topseudosquare<0.8cm,0.8cm,#6>$}
\vskip-1.4pt
\hbox{
        $\Tophalfsquare<0.8cm,0.8cm,#7>
         \topseudosquare(0.8cm,0.8cm,#8)
         \topseudosquare(0.8cm,0.8cm,#7)
         \Topseudosquare<0.8cm,0.8cm,#7>$}
}}

\def\sevenfourA(#1,#2,#3,#4,#5,#6,#7,#8,#9){
\vcenter{\hbox{
        $\tophalfsquare(0.8cm,1.5cm,#1)
         \topseudosquare(0.8cm,1.5cm,#8)
         \topseudosquare(0.8cm,1.5cm,#2)$}
\vskip-1.4pt
\hbox{
        $\tophalfsquare(1.0cm,1.5cm,#9)
         \topseudosquare(1.0cm,1.5cm,#9)
         \topseudosquare(1.0cm,1.5cm,#9)$}
\vskip-1.2pt
\hbox{
        $\tophalfsquare(0.8cm,1.5cm,#3)
         \topseudosquare(0.8cm,1.5cm,#8)
         \topseudosquare(0.8cm,1.5cm,#4)$}
\vskip-1.2pt
\hbox{
        $\tophalfsquare(0.8cm,1.5cm,#5)
         \topseudosquare(0.8cm,1.5cm,#8)
         \topseudosquare(0.8cm,1.5cm,#5)$}
\vskip-1.2pt
\hbox{
        $\tophalfsquare(0.8cm,1.5cm,#6)
         \topseudosquare(0.8cm,1.5cm,#8)
         \topseudosquare(0.8cm,1.5cm,#6)$}
\vskip-1.2pt
\hbox{
        $\tophalfsquare(1.0cm,1.5cm,#9)
         \topseudosquare(1.0cm,1.5cm,#9)
         \topseudosquare(1.0cm,1.5cm,#9)$}
\vskip-1.2pt
\hbox{
        $\tophalfsquare(0.8cm,1.5cm,#7)
         \topseudosquare(0.8cm,1.5cm,#8)
         \topseudosquare(0.8cm,1.5cm,#7)$}
}}

\def\sevenfourB(#1,#2,#3,#4,#5,#6,#7,#8,#9){
\vcenter{\hbox{
        $\tophalfsquare(0.8cm,3.0cm,#1)
         \topseudosquare(0.8cm,3.0cm,#8)
         \topseudosquare(0.8cm,3.0cm,#2)$}
\vskip-1.4pt
\hbox{
        $\tophalfsquare(1.0cm,3.0cm,#9)
         \topseudosquare(1.0cm,3.0cm,#9)
         \topseudosquare(1.0cm,3.0cm,#9)$}
\vskip-1.2pt
\hbox{
        $\tophalfsquare(0.8cm,3.0cm,#3)
         \topseudosquare(0.8cm,3.0cm,#8)
         \topseudosquare(0.8cm,3.0cm,#4)$}
\vskip-1.2pt
\hbox{
        $\tophalfsquare(0.8cm,3.0cm,#5)
         \topseudosquare(0.8cm,3.0cm,#8)
         \topseudosquare(0.8cm,3.0cm,#5)$}
\vskip-1.2pt
\hbox{
        $\tophalfsquare(0.8cm,3.0cm,#6)
         \topseudosquare(0.8cm,3.0cm,#8)
         \topseudosquare(0.8cm,3.0cm,#6)$}
\vskip-1.2pt
\hbox{
        $\tophalfsquare(1.0cm,3.0cm,#9)
         \topseudosquare(1.0cm,3.0cm,#9)
         \topseudosquare(1.0cm,3.0cm,#9)$}
\vskip-1.2pt
\hbox{
        $\tophalfsquare(0.8cm,3.0cm,#7)
         \topseudosquare(0.8cm,3.0cm,#8)
         \topseudosquare(0.8cm,3.0cm,#7)$}
}}

\def\seventhreeA(#1,#2,#3,#4,#5,#6,#7,#8,#9){
\vcenter{\hbox{
        $\tophalfsquare(0.8cm,1.6cm,#1)
         \topseudosquare(0.8cm,1.6cm,#8)
         \topseudosquare(0.8cm,1.6cm,#2)$}
\vskip-1.4pt
\hbox{
        $\tophalfsquare(1.0cm,1.6cm,#9)
         \topseudosquare(1.0cm,1.6cm,#9)
         \topseudosquare(1.0cm,1.6cm,#9)$}
\vskip-1.2pt
\hbox{
        $\tophalfsquare(0.8cm,1.6cm,#3)
         \topseudosquare(0.8cm,1.6cm,#8)
         \topseudosquare(0.8cm,1.6cm,#4)$}
\vskip-1.2pt
\hbox{
        $\tophalfsquare(0.8cm,1.6cm,#5)
         \topseudosquare(0.8cm,1.6cm,#8)
         \topseudosquare(0.8cm,1.6cm,#5)$}
\vskip-1.2pt
\hbox{
        $\tophalfsquare(0.8cm,1.6cm,#6)
         \topseudosquare(0.8cm,1.6cm,#8)
         \topseudosquare(0.8cm,1.6cm,#6)$}
\vskip-1.2pt
\hbox{
        $\tophalfsquare(1.0cm,1.6cm,#9)
         \topseudosquare(1.0cm,1.6cm,#9)
         \topseudosquare(1.0cm,1.6cm,#9)$}
\vskip-1.2pt
\hbox{
        $\tophalfsquare(0.8cm,1.6cm,#7)
         \topseudosquare(0.8cm,1.6cm,#8)
         \topseudosquare(0.8cm,1.6cm,#7)$}
}}

\title{\textbf{\large{Explicit Forms of Cluster Variables 
on Double Bruhat Cells $G^{u,e}$ of type C}}}
\author{\normalsize{YUKI KANAKUBO\thanks{Division of Mathematics, 
Sophia University, Kioicho 7-1, Chiyoda-ku, Tokyo 102-8554,
Japan: {j\_chi\_sen\_you\_ky@eagle.sophia.ac.jp}}
\ and\ 
TOSHIKI NAKASHIMA\thanks{Division of Mathematics, 
Sophia University, Kioicho 7-1, Chiyoda-ku, Tokyo 102-8554,
Japan: { toshiki@sophia.ac.jp}:
supported in part by JSPS Grants 
in Aid for Scientific Research $\sharp 22540031$.}
}}
\date{}

\begin{document}

\maketitle
\vspace{-10pt}
\centerline{{\it Dedicated to Professor Ken-ichi SHINODA}}

\begin{abstract}
Let $G=Sp_{2r}({\mathbb C})$
be a simply connected simple algebraic group over $\mathbb{C}$ of type $C_r$, 
$B$ and $B_-$ its two opposite Borel subgroups, and $W$ the associated 
 Weyl group. 
For $u$, $v\in W$, 
it is known that the coordinate ring 
${\mathbb C}[G^{u,v}]$ of the double Bruhat cell $G^{u,v}=BuB\cap B_-vB_-$ is 
isomorphic to an upper cluster algebra $\ovl{\cA}(\textbf{i})_{{\mathbb C}}$ 
and the generalized minors $\Delta(k;\textbf{i})$ are 
the cluster variables of ${\mathbb C}[G^{u,v}]$\cite{A-F-Z}. 
In the case $v=e$, 
we shall describe the generalized minor $\Delta(k;\textbf{i})$
 explicitly.
\end{abstract}


\tableofcontents

\section{Introduction}

Let $G$ be a simply connected simple algebraic group over $\bbC$ of rank
$r$, $B,B_-\subset G$ the opposite Borel subgroups, $H:=B\cap B^-$ the maximal torus, 
$N\subset B$, $N_-\subset B_-$ the maximal unipotent 
subgroups and $W$ the associated Weyl group. For $u,v\in W$, define
$G^{u,v}:=(Bu B)\cap(B_- v B_-)$ (resp. $L^{u,v}:=(Nu N)\cap(B_- v B_-)$)
and call it the (reduced) double Bruhat cell.

In \cite{A-F-Z}, it is shown that 
 the coordinate ring $\bbC[G^{u,v}]$ ($u,v\in W$) of double Bruhat
cell $G^{u,v}$ has the structure of an upper cluster algebra.
The initial cluster variables of this upper cluster algebras are given
as certain generalized minors on $G^{u,v}$.

In \cite{KaN}, we treated the case of type A and $v=e$, where we described 
the explicit forms of the generalized minors $\{\Delta(k;\textbf{i})\}$ and 
revealed the linkage between $\Delta(k;\textbf{i})$ and monomial realizations 
of crystals.

In this paper, we shall write down the explicit forms of the generalized minors
$\Delta(k;\textbf{i})$
on the (reduced) double Bruhat cell $G^{u,e}$ ($L^{u,e}$) of type $C_r$
by using the
`path descriptions' of generalized minors (see Sect.6), where 
we only treat a Weyl group elements $u$ with the form as in 
\eqref{uset0} and denote its reduced word $\bfi$ by \eqref{iset0}.
Indeed, generalized minors are expressed in terms of certain invariant
 bilinear forms (see \eqref{minor-bilin}). And then, using this bilinearity 
we obtain `path descriptions' of the generalized minors.

Unfortunately, we do not present the relation between the explicit forms of 
$\Delta(k;\textbf{i})$ and crystals here unlike with \cite{KaN}.
We will, however, discuss this elsewhere.

The main result is given as in Theorem \ref{thm1}:
Let $\bfi$ be the reduced word of $u$ as above and $i_k$ is the $k$-th 
index of $\bfi$ form the left. In \cite{B-Z}, it is shown that there exists a biregular isomorphism from
$(\mathbb{C}^{\times})^{n}$ to a Zariski open subset of $L^{u,e}$ $(n:=l(u))$ (see Theorem \ref{fp2}). We denote this isomorphism by $x^L_{\textbf{i}}$ and set $\Delta^L(k;\textbf{i}):=\Delta(k;\textbf{i})\circ x^L_{\textbf{i}}$. We also set the monomials $\ovl C(l,k)$ and $C(l,k)$ as in \eqref{ccbar}.
\\
{\bf Theorem 5.7} 
We set $d:=i_k=i_n$ and
\[ \textbf{Y}:=(Y_{1,1},Y_{1,2},\cdots,Y_{1,r},\cdots,Y_{m-1,1},\cdots,Y_{m-1,r},Y_{m,1},\cdots,Y_{m,i_n})\in  (\mathbb{C}^{\times})^{n}. \]
Then we have
\begin{multline*} 
\Delta^L(k;\textbf{i})(\textbf{Y})= \sum_{(*)}\prod^{d}_{i=1} \ovl{C}(m-l^{(1)}_i,k^{(1)}_i)\cdots \ovl{C}(m-l^{(\delta_i)}_i,k^{(\delta_i)}_i)\\
\cdot C(m-l^{(\delta_i+1)}_i,|k^{(\delta_i+1)}_i|-1)\cdots C(m-l^{(m-m')}_i,|k^{(m-m')}_i|-1),
\end{multline*} 
\[ l^{(s)}_{i}:=
\begin{cases}
k^{(s)}_{i}+s-i-1 & {\rm if}\  s\leq \del_{i}, \\
s-i+r & {\rm if}\  s>\del_{i}
\end{cases}
\q (1\leq i\leq d)
\]
where $(*)$ is the conditions for 
$k^{(s)}_i$ $(1\leq s\leq m-m',\ 1\leq i\leq d)$ : 
$1\leq k^{(s)}_1<k^{(s)}_2<\cdots<k^{(s)}_d\leq\ovl{1}\q (1\leq s\leq m-m')$,
$1\leq k^{(1)}_i\leq \cdots\leq k^{(m-m')}_i\leq m'+i\q (1\leq i\leq r-m')$,
and $1\leq k^{(1)}_i\leq \cdots\leq k^{(m-m')}_i\leq \ovl{1}\q 
(r-m'+1\leq i\leq d)$, and 
$\delta_i$ $(i=1,\cdots,d)$ are the numbers which satisfy 
$1\leq k^{(1)}_i\leq k^{(2)}_i\leq \cdots\leq k^{(\delta_i)}_i\leq r$, 
$\ovl{r}\leq k^{(\delta_i+1)}_i\leq \cdots\leq k^{(m-m')}_i\leq\ovl{1}$.

For ${\bf k}=(k^{(s)}_i)$ and ${\bf k}'=(k^{'(s)}_i)$ satisfying ($*$), let us 
write the monomial
\begin{eqnarray*}
&&C({\bf k}):=\prod^{d}_{i=1} \ovl{C}(m-l^{(1)}_i,k^{(1)}_i)\cdots \ovl{C}(m-l^{(\delta_i)}_i,k^{(\delta_i)}_i)\\
&&\qq \cdot C(m-l^{(\delta_i+1)}_i,|k^{(\delta_i+1)}_i|-1)\cdots C(m-l^{(m-m')}_i,|k^{(m-m')}_i|-1).
\end{eqnarray*}
Note that even if ${\bf k}\ne{\bf k}'$, we may have $C(\bf k)=C(\bf k')$.
Thus, we will know that the coefficients of the monomials in $\Del^L(k;{\bf i})$
are not necessarily 1 (See Example \ref{pathex3}).
We shall show Theorem \ref{thm1} in the last section by using 
``path descriptions''. 
By this theorem, we find that all the generalized minors 
$\{\Delta^L(k;\textbf{i})(\textbf{Y})\}$ are Laurent polynomials with 
non-negative coefficients.

Finally, we also define $\Delta^G(k;\textbf{i}):=\Delta(k;\textbf{i})\circ \overline{x}^G_{\textbf{i}}$, where $\overline{x}^G_{\textbf{i}}$ is a biregular isomorphism from
$H\times(\mathbb{C}^{\times})^{n}$ to a Zariski open subset of $G^{u,e}$ (see Proposition \ref{gprime}). In Proposition \ref{gprop}, we shall show that $\Delta^G(k;\textbf{i})$ is immediately obtained from $\Delta^L(k;\textbf{i})$.


\section{Fundamental representations for type $C_r$}

Let $I:=\{1,\cdots,r\}$ be a finite index set and $A=(a_{ij})_{i,j\in I}$ 
be the Cartan matrix of type $C_r$.
That is, $A=(a_{i,j})_{i,j\in I}$ is given by
\[a_{i,j}=
\begin{cases}
2 & {\rm if}\ i=j, \\
-1 & {\rm if}\ |i-j|=1\ {\rm and}\ (i,j)\neq (r-1,r), \\
-2 & {\rm if}\ (i,j)=(r-1,r), \\
0 & {\rm otherwise.}  
\end{cases}
\]
Let $(\frh,\{\al_i\}_{i\in I},\{h_i\}_{i\in I})$ 
be the associated
root data 
satisfying $\al_j(h_i)=a_{ij}$ where 
$\al_i\in \frh^*$ is a simple root and 
$h_i\in \frh$ is a simple co-root.
Note that $\al_i\ (i\neq r)$ are short roots and $\al_r$ is the long root. 
Let $\{\Lm_i\}_{i\in I}$ be the set of the fundamental 
weights satisfying $\al_j(h_i)=a_{i,j}$ 
and $\Lm_i(h_j)=\del_{i,j}$. Let $P=\bigoplus_{i\in I}\mathbb{Z}\Lm_i$ be the weight lattice and $P^*=\bigoplus_{i\in I}\mathbb{Z}h_i$ be the dual weight lattice. Define the order on the set $J:=\{i,\ovl i|1\leq i\leq r\}$ by 
\begin{equation}\label{order}
 1< 2<\cd< r-1< r
< \ovl r< \ovl{r-1}< \cd< \ovl
 2< \ovl 1.
\end{equation}
For the simple Lie algebra 
 $\frg=\frs\frp(2r,\mathbb{C})=\lan \frh,e_i,f_i(i\in I)\ran$, 
let us describe the vector representation 
$V(\Lm_1)$. Set ${\mathbf B}^{(r)}:=
\{v_i,v_{\ovl i}|i=1,2,\cd,r\}$ and define 
$V(\Lm_1):=\bigoplus_{v\in{\mathbf B}^{(r)}}\bbC v$. The weights of $v_i$, $v_{\ovl{i}}$ $(i=1,\cd,r)$ are as follows:
\begin{equation}\label{wtv}
 {\rm wt}(v_i)=\Lm_i-\Lm_{i-1},\q {\rm wt}(v_{\ovl{i}})=\Lm_{i-1}-\Lm_{i},
\end{equation}
where $\Lm_0=0$. We define the $\frs\frp(2r,\mathbb{C})$-action on $V(\Lm_1)$ as follows:
\begin{eqnarray}
&& h v_j=\lan h,{\rm wt}(v_j)\ran v_j\ \ (h\in P^*,\ j\in J), \\
&&f_iv_i=v_{i+1},\ f_iv_{\ovl{i+1}}=v_{\ovl i},\q
e_iv_{i+1}=v_i,\ e_iv_{\ovl i}=v_{\ovl{i+1}}
\q(1\leq i<r),\label{c-f1}\\
&&f_r v_r=v_{\ovl r},\qq 
e_r v_{\ovl r}=v_r,\label{c-f2}
\end{eqnarray}
and the other actions are trivial.

Let $\Lm_i$ be the $i$-th fundamental weight of type $C_r$.
As is well-known that the fundamental representation 
$V(\Lm_i)$ $(1\leq i\leq r)$
is embedded in $\wedge^i V(\Lm_1)$
with multiplicity free.
The explicit form of the highest (resp. lowest) weight 
vector $u_{\Lm_i}$ (resp. $v_{\Lm_i}$)
of $V(\Lm_i)$ is realized in 
$\wedge^i V(\Lm_1)$ as follows:
\begin{equation}
\begin{array}{ccc}\displaystyle
u_{\Lm_i}&=&v_1\wedge v_2\wedge\cdots\wedge v_i,\\
v_{\Lm_i}&=&v_{\ovl{1}}\wedge v_{\ovl{2}}\wedge\cdots \wedge v_{\ovl{i}}.
\end{array}
\label{h-l}
\end{equation}

\section{Factorization theorem for type C}\label{DBCs}

In this section, we shall introduce (reduced) double Bruhat cells $G^{u,v}$, $L^{u,v}$, and their properties in the case $G=Sp_{2r}({\mathbb C})$, $v=e$ and some special $u\in W$. In \cite{B-Z} and \cite{F-Z}, these properties had been proven for simply connected, connected, semisimple complex algebraic groups and arbitrary $u,v\in W$. 

For $l\in \mathbb{Z}_{>0}$, we set $[1,l]:=\{1,2,\cdots,l\}$.

\subsection{Double Bruhat cells}\label{factpro}

Let $G=Sp_{2r}(\mathbb{C})$ be the simple complex algebraic
group of type ${\rm C}_r$, $B$ and $B_-$ two opposite Borel subgroups in $G$, $N\subset B$ and $N_-\subset B_-$ their unipotent radicals, 
$H:=B\cap B_-$ a maximal torus. We set $\frg:={\rm Lie}(G)$ with the Cartan decomposition $\frg=\frn_-\oplus \frh \oplus \frn$. Let $e_i$, $f_i$ $(i\in[1,r])$ be the generators of $\frn$, $\frn_-$. For $i\in[1,r]$ and $t \in \mathbb{C}$, we set $x_i(t):={\rm exp}(te_i),\ y_{i}:={\rm exp}(tf_i)$. Let $\varphi_i :SL_2(\mathbb{C}) \rightarrow G$ be the canonical embedding corresponding to simple root $\alpha_i$. Then we have
\begin{equation}\label{xiyidef} 
x_i(t)=\varphi_i 
\begin{pmatrix}
1 & t \\
0 & 1 
\end{pmatrix},\q
y_{i}(t)=\varphi_i 
\begin{pmatrix}
1 & 0 \\
t & 1 
\end{pmatrix}.
\end{equation}
Let $W:=\lan s_i |i=1,\cdots,r \ran$ be the Weyl group of $\frg$, where
$\{s_i\}$ are the simple reflections. We identify the Weyl group $W$ with ${\rm Norm}_G(H)/H$. An element 
\begin{equation}\label{smpl}
\ovl{s_i}:=x_i(-1)y_i(1)x_i(-1)
\end{equation}
is in ${\rm Norm}_G(H)$, which is representative of $s_i\in W={\rm Norm}_G(H)/H$ \cite{N1}. For $u\in W$, we denote the length of $u$ by $l(u)$.

We have two kinds of Bruhat decompositions of $G$ as follows:
\[ G=\displaystyle\coprod_{u \in W}BuB=\displaystyle\coprod_{u \in W}B_-uB_- .\]
Then, for $u$, $v\in W$, 
we define the {\it double Bruhat cell} $G^{u,v}$ as follows:
\[ G^{u,v}:=BuB \cap B_-vB_-. \]
This is biregularly isomorphic to a Zariski open subset of 
an affine space of dimension $r+l(u)+l(v)$ \cite[Theorem 1.1]{F-Z}.

We also define the {\it reduced double Bruhat cell} $L^{u,v}$ as follows:
\[ L^{u,v}:=NuN \cap B_-vB_- \subset G^{u,v}. \] 
As is similar to the case $G^{u,v}$, $L^{u,v}$ is 
biregularly isomorphic to a Zariski open subset of an 
affine space of dimension $l(u)+l(v)$ \cite[Proposition 4.4]{B-Z}.

\begin{defn}\label{redworddef}
Let $u=s_{i_1}\cdots s_{i_n}$ be a reduced expression of $u\in W$ $(i_1,\cdots,i_n\in [1,r])$. Then the finite sequence 
\[ \textbf{i}:=(i_1,\cdots,i_n) \]
is called a {\it reduced word} for $u$.
\end{defn}

In this paper, we treat (reduced) Double Bruhat cells of the form $G^{u,e}:=BuB \cap B_-$ and $L^{u,e}:=NuN \cap B_-$, where $u\in W$ is an element whose reduced word can be written as a left factor of $(1,2,3,\cdots,r)^r$:
\begin{equation}\label{uset0}
u=(s_1s_2\cdots s_r)^{m-1}s_1\cdots s_{i_n},
\end{equation} 
where $n:=l(u)$ is the length of $u$ and $1\leq i_n\leq r$. Let $\textbf{i}$ be a reduced word of $u$:
\begin{equation}\label{iset0}
 \textbf{i}=(\underbrace{1,\cdots,r}_{1{\rm \,st\
cycle}},\underbrace{1,\cdots,r}_{2{\rm \,nd\
cycle}},\cdots,\underbrace{1,\cdots,r}_{m-1{\rm \,\,th\ cycle}},
\underbrace{1,2,\cdots,i_n}_{m{\rm \,\,th\ cycle}}).
\end{equation}
Note that $(1,2,3,\cdots,r)^r$ is a reduced word of the longest element in $W$.

\subsection{Factorization theorem for type $C_r$}\label{factproC}

In this subsection, we shall introduce the isomorphisms between double Bruhat cell $G^{u,e}$ and $H\times (\mathbb{C}^{\times})^{l(u)}$, and between $L^{u,e}$ and $(\mathbb{C}^{\times})^{l(u)}$. 
As in the previous section, we consider the case $G=Sp_{2r}(\mathbb{C})$. 

For a reduced word $\textbf{i}=(i_1, \cdots ,i_n)$ 
($i_1,\cdots,i_n\in[1,r]$), 
we define a map $x^G_{\textbf{i}}:H\times \mathbb{C}^n \rightarrow G$ as 
\begin{equation}\label{xgdef}
x^G_{\textbf{i}}(a; t_1, \cdots, t_n):=a\cdot y_{i_1}(t_1)\cdots y_{i_n}(t_n).
\end{equation}

\begin{thm}\label{fp}${\cite[Theorem\ 1.2]{F-Z}}$ We set $u\in W$ and its reduced word $\textbf{i}$ as in $(\ref{uset0})$ and $(\ref{iset0})$. The map $x^G_{\textbf{i}}$ defined above can be restricted to a biregular isomorphism between $H\times (\mathbb{C}^{\times})^{l(u)}$ and a Zariski open subset of $G^{u,e}$. 
\end{thm}

Next, for $i \in [1,r]$ and $t\in \mathbb{C}^{\times}$, we define as follows:
\begin{equation}\label{alxmdef}
\alpha_i^{\vee}(t):=t^{h_i}=\varphi_i 
\begin{pmatrix}
t & 0 \\
0 & t^{-1} 
\end{pmatrix},\ \ 
 x_{-i}(t):=y_{i}(t)\alpha_i^{\vee}(t^{-1})=\varphi_i 
\begin{pmatrix}
t^{-1} & 0 \\
1 & t
\end{pmatrix}.
\end{equation}
For $\textbf{i}=(i_1, \cdots ,i_n)$
($i_1,\cdots,i_n\in[1,r]$), 
we define a map $x^L_{\textbf{i}}:\mathbb{C}^n \rightarrow G$ as 
\begin{equation}\label{xldef}
x^L_{\textbf{i}}(t_1, \cdots, t_n):=x_{-i_1}(t_1)\cdots x_{-i_n}(t_n).
\end{equation}
We have the following theorem which is similar to the previous one.
\begin{thm}\label{fp2}${\cite[Proposition\ 4.5]{B-Z}}$
We set $u\in W$ and its reduced word $\textbf{i}$ as in $(\ref{uset0})$ and $(\ref{iset0})$.
The map $x^L_{\textbf{i}}$ defined above can be restricted to a
biregular isomorphism between $ (\mathbb{C}^{\times})^{l(u)}$ 
and a Zariski open subset of $L^{u,e}$. 
\end{thm}

We define a map
$\ovl{x}^G_{\textbf{i}}:H\times(\mathbb{C}^{\times})^{n}\rightarrow
G^{u,e}$ as
\[ \ovl{x}^G_{\textbf{i}}(a;t_1,\cdots,t_n)
=ax^L_{\textbf{i}}(t_1,\cdots,t_n), \]
where $a\in H$ and $(t_1,\cdots,t_n)\in (\mathbb{C}^{\times})^{n}$.
\begin{prop}\label{gprime}
In the above setting, the map $\ovl{x}^G_{\textbf{i}}$ is a biregular isomorphism between $H\times(\mathbb{C}^{\times})^{n}$ and a Zariski open subset of $G^{u,e}$.
\end{prop}

\nd
{\sl Proof.}
In this proof, we use the notation
\[ (Y_{1,1},\cdots,Y_{1,r},\cdots,Y_{m-1,1},\cdots,Y_{m-1,r},Y_{m,1},\cdots,Y_{m,i_n})\in (\mathbb{C}^{\times})^{n} \]
for variables instead of $(t_1,\cdots,t_n)$.

We define a map
$\phi:H\times(\mathbb{C}^{\times})^{n}\rightarrow
H\times(\mathbb{C}^{\times})^{n}$ as follows: For 
\[
\textbf{Y}:=(a;Y_{1,1},\cdots,Y_{1,r}, \cdots,Y_{m,1},\cdots,Y_{m,i_n}),
\]
we define
$\phi(\textbf{Y})=(\Phi_a(\textbf{Y});\Phi_{1,1}(\textbf{Y}),\cdots,\Phi_{1,r}(\textbf{Y}),\cdots,\Phi_{m,1}(\textbf{Y}),\cdots,\Phi_{m,i_n}(\textbf{Y}))$ as
\[ \Phi_a(\textbf{Y}):=a\cdot\left(\prod^{m-1}_{j=1}{\al_1^{\vee}(Y_{j,1})^{-1}\cdots
\al_r^{\vee}(Y_{j,r})^{-1}}\right)\cdot\al_1^{\vee}(Y_{m,1})^{-1}\cdots
\al_{i_n}^{\vee}(Y_{m,i_n})^{-1}, \]
and for $1\leq s\leq m$,
\begin{equation}\label{mbase0} 
\Phi_{s,l}(\textbf{Y}):=
\begin{cases}
\frac{(Y_{s+1,l-1}Y_{s+2,l-1}\cdots Y_{m,l-1})(Y_{s,l+1}Y_{s+1,l+1}\cdots Y_{m,l+1})}{Y_{s,l}(Y_{s+1,l}\cdots Y_{m,l})^{2}} & {\rm if}\ 1\leq l<r, \\
\frac{(Y_{s+1,r-1}Y_{s+2,r-1}\cdots Y_{m,r-1})^2}{Y_{s,r}(Y_{s+1,r}\cdots Y_{m,r})^{2}} & {\rm if}\ l=r,
\end{cases}
\end{equation}
where in (\ref{mbase0}), if we see the variables $Y_{\zeta,0}$ $(1\leq\zeta\leq m)$ and $Y_{m,\xi}$ $(i_n<\xi)$, then we understand $Y_{\zeta,0}=Y_{m,\xi}=1$. For example, $Y_{s+1,l-1}=1$ in the case $l=1$.
Note that $\phi$ is a biregular isomorphism since we can recurrently construct the inverse map $\psi:H\times(\mathbb{C}^{\times})^{n}\rightarrow
H\times(\mathbb{C}^{\times})^{n}$, $\textbf{Y}\mapsto (\Psi_a(\textbf{Y});\Psi_{1,1}(\textbf{Y}),\cdots,\Psi_{m,i_n}(\textbf{Y}))$ of $\phi$ as follows: The definition (\ref{mbase0}) implies that $\Phi_{m,i_n}(\textbf{Y})=\frac{1}{Y_{m,i_n}}$, and hence $Y_{m,i_n}=\frac{1}{\Psi_{m,i_n}(\textbf{Y})}$. So we set $\Psi_{m,i_n}(\textbf{Y})=\frac{1}{Y_{m,i_n}}$. Suppose that we can construct $\Psi_{m,i_n}(\textbf{Y}),\Psi_{m,i_n-1}(\textbf{Y}),\cdots$,
$\Psi_{m,1}(\textbf{Y}),\cdots\Psi_{s+1,r}(\textbf{Y})$,
$\cdots,\Psi_{s+1,1}(\textbf{Y}),\Psi_{s,r}(\textbf{Y}),\cdots$,
$\Psi_{s,l+1}(\textbf{Y})$. Then we define
\[ \Psi_{s,l}(\textbf{Y}):=
\begin{cases}
\frac{(\Psi_{s+1,l}(\textbf{Y})\cdots \Psi_{m,l}(\textbf{Y}))^{2}}
{Y_{s,l}(\Psi_{s+1,l-1}(\textbf{Y})\Psi_{s+2,l-1}(\textbf{Y})\cdots \Psi_{m,l-1}(\textbf{Y}))(\Psi_{s,l+1}(\textbf{Y})\cdots \Psi_{m,l+1}(\textbf{Y}))}
 & {\rm if}\ 1\leq l<r, \\
\frac{(\Psi_{s+1,r}(\textbf{Y})\cdots \Psi_{m,r}(\textbf{Y}))^{2}}
{Y_{s,r}(\Psi_{s+1,r-1}(\textbf{Y})\Psi_{s+2,r-1}(\textbf{Y})\cdots \Psi_{m,r-1}(\textbf{Y}))^2}
 & {\rm if}\ l=r.
\end{cases}
\]
We also define
\[ \Psi_{a}(\textbf{Y}):=
a\cdot\left(\prod^{m-1}_{j=1}{\al_1^{\vee}(\Psi_{j,1}(\textbf{Y}))\cdots
\al_r^{\vee}(\Psi_{j,r}(\textbf{Y}))}\right)\cdot\al_1^{\vee}(\Psi_{m,1}(\textbf{Y}))\cdots
\al_{i_n}^{\vee}(\Psi_{m,i_n}(\textbf{Y})).
\]
Then, we get the inverse map $\psi$ of $\phi$.

Let us prove
\[ \ovl{x}^G_{\textbf{i}}(\textbf{Y})=(x^G_{\textbf{i}}\circ\phi)(\textbf{Y}), \]
which implies that $\ovl{x}^G_{\textbf{i}}:H\times(\mathbb{C}^{\times})^{n}\rightarrow G^{u,e}$ is biregular isomorphism by Theorem \ref{fp}.

First, it is known that
\begin{equation}\label{base2}
\al_i^{\vee}(c)^{-1}y_{j}(t)=\begin{cases}
	y_{i}(c^2t)\al_i^{\vee}(c)^{-1} & {\rm if}\ i=j, \\
	y_{j}(c^{-1}t)\al_i^{\vee}(c)^{-1} & {\rm if}\ |i-j|=1\ {\rm and}\ (i,j)\neq (r-1,r), \\
	y_{j}(c^{-2}t)\al_i^{\vee}(c)^{-1} & {\rm if}\ (i,j)= (r-1,r), \\
	y_{j}(t)\al_i^{\vee}(c)^{-1} & {\rm otherwise},
\end{cases}
\end{equation}
for $1\leq i,\ j\leq r$ and $c,\ t\in \mathbb{C}^{\times}$.

On the other hand, it follows from the definition (\ref{xgdef}) of $x^G_{\textbf{i}}$ and $(\ref{mbase0})$ that
\begin{multline}\label{xigp}(x^G_{\textbf{i}}\circ\phi)(\textbf{Y})\\
=a\times\left(\prod^{m-1}_{j=1}{\al_1^{\vee}(Y_{j,1})^{-1}\cdots
\al_r^{\vee}(Y_{j,r})^{-1}}\right)\cdot\al_1^{\vee}(Y_{m,1})^{-1}\cdots
\al_{i_n}^{\vee}(Y_{m,i_n})^{-1}
\\ 
 \times y_{1}(\Phi_{1,1}(\textbf{Y}))y_{2}(\Phi_{1,2}(\textbf{Y}))\cdots
 y_{r}(\Phi_{1,r}(\textbf{Y}))\cdots y_{1}(\Phi_{m,1}(\textbf{Y}))
\cdots y_{i_n}(\Phi_{m,i_n}(\textbf{Y})).
\end{multline}

For each $s$ and $l$ $(1\leq s\leq m,\ 1\leq l\leq r)$, we can move 
\begin{multline*}
\al_{l}^{\vee}(Y_{s,l})^{-1}
\al_{l+1}^{\vee}(Y_{s,l+1})^{-1}
\cdots\al_{r}^{\vee}(Y_{s,r})^{-1}\\
\cdot\left(\prod^{m-1}_{j=s+1}\al_1^{\vee}(Y_{j,1})^{-1}\cdots \al_r^{\vee}(Y_{j,r})^{-1}\right)\cdot
\al_1^{\vee}(Y_{m,1})^{-1}\cdots\al_{i_n}^{\vee}(Y_{m,i_n})^{-1}
\end{multline*}
to the right of
$y_{l}(\Phi_{s,l}(\textbf{Y}))$ by using the relations (\ref{base2}). 
For example,
\begin{flushleft}
$\al_1^{\vee}(Y_{m,1})^{-1}\cdots
\al_{i_n}^{\vee}(Y_{m,i_n})^{-1} y_{l}(\Phi_{s,l}(\textbf{Y}))=$
\end{flushleft}
\[
\begin{cases}
y_{l}\left(\frac{Y_{m,l}^2}{Y_{m,l-1}Y_{m,l+1}}\Phi_{s,l}(\textbf{Y})\right)\al_1^{\vee}(Y_{m,1})^{-1}\cdots
\al_{i_n}^{\vee}(Y_{m,i_n})^{-1} &\ {\rm if}\ 1\leq l<r, \\
y_{r}\left(\frac{Y_{m,r}^2}{Y_{m,r-1}^2}\Phi_{s,r}(\textbf{Y})\right)\al_1^{\vee}(Y_{m,1})^{-1}\cdots
\al_{i_n}^{\vee}(Y_{m,i_n})^{-1} &\ {\rm if}\ l=r.
\end{cases}
\]

Repeating this argument, in the case $l<r$, we have
\begin{multline*}
\al_{l}^{\vee}(Y_{s,l})^{-1}
\al_{l+1}^{\vee}(Y_{s,l+1})^{-1}
\cdots\al_{r}^{\vee}(Y_{s,r})^{-1}\\
\times\left(\prod^{m-1}_{j=s+1}\al_1^{\vee}(Y_{j,1})^{-1}\cdots \al_r^{\vee}(Y_{j,r})^{-1}\right)\cdot
\al_1^{\vee}(Y_{m,1})^{-1}\cdots\al_{i_n}^{\vee}(Y_{m,i_n})^{-1}y_{l}(\Phi_{s,l}(\textbf{Y}))\\
=y_{l}\left(\frac{(Y_{s,l}Y_{s+1,l}\cdots Y_{m-1,l}Y_{m,l})^2}{(Y_{s+1,l-1}\cdots Y_{m-1,l-1}Y_{m,l-1})(Y_{s,l+1}\cdots Y_{m-1,l+1}Y_{m,l+1})}\Phi_{s,l}(\textbf{Y})\right)\cdot \al_{l}^{\vee}(Y_{s,l})^{-1}\\
\times\al_{l+1}^{\vee}(Y_{s,l+1})^{-1}
\cdots\al_{r}^{\vee}(Y_{s,r})^{-1}
\cdot\left(\prod^{m-1}_{j=s+1}\al_1^{\vee}(Y_{j,1})^{-1}\cdots \al_r^{\vee}(Y_{j,r})^{-1}\right)\cdot
\al_1^{\vee}(Y_{m,1})^{-1}\cdots\al_{i_n}^{\vee}(Y_{m,i_n})^{-1}.
\end{multline*}

Note that $\frac{(Y_{s,l}Y_{s+1,l}\cdots Y_{m-1,l}Y_{m,l})^2}{(Y_{s+1,l-1}\cdots Y_{m-1,l-1}Y_{m,l-1})(Y_{s,l+1}\cdots Y_{m-1,l+1}Y_{m,l+1})}\Phi_{s,l}(\textbf{Y})=Y_{s,l}$, which implies 
\begin{equation}\label{tradeeq}
\al_{l+1}^{\vee}(Y_{s,l+1})^{-1}\cdots \al_{i_n}^{\vee}(Y_{m,i_n})^{-1}y_{l}(\Phi_{s,l}(\textbf{Y}))=y_{l}(Y_{s,l})\al_{l+1}^{\vee}(Y_{s,l+1})^{-1}\cdots \al_{i_n}^{\vee}(Y_{m,i_n})^{-1}.
\end{equation}

In the case $l=r$, we can also verify the relation (\ref{tradeeq}) similarly. Thus, by (\ref{xigp}) and (\ref{tradeeq}), we have
\begin{eqnarray*}
&&(x^G_{\textbf{i}}\circ\phi)(\textbf{Y})=a\cdot
y_{1}(Y_{1,1})\al_1^{\vee}(Y_{1,1})^{-1}\cdots y_{r}(Y_{1,r})\al_r^{\vee}(Y_{1,r})^{-1}\times\cdots \\
&&\qq\times y_{1}(Y_{m,1})\al_1^{\vee}(Y_{m,1})^{-1}\cdots
 y_{i_n}(Y_{m,i_n})\al_{i_n}^{\vee}(Y_{m,i_n})^{-1}\\
&&
=a\cdot x_{-1}(Y_{1,1})\cdots x_{-r}(Y_{1,r})\cdots x_{-1}(Y_{m,1})\cdots
x_{-i_n}(Y_{m,i_n})\\
&&=\ovl{x}^G_{\textbf{i}}(\textbf{Y}).\qq\qq\qq\qq\qq\qq\qq\qq\qq \qed
\end{eqnarray*}

\section{Cluster algebras and generalized minors}
For this section, see {\it e.g.,}\cite{M-M-A,F-Z,FZ2,A-F-Z}.

We set $[1,l]:=\{1,2,\cdots,l\}$ and $[-1,-l]:=\{-1,-2,\cdots,-l\}$ for $l\in \mathbb{Z}_{>0}$. For $n,m\in \mathbb{Z}_{>0}$, let $x_1, \cdots ,x_n,x_{n+1}, \cdots
,x_{n+m}$ be variables and $\mathcal{P}$ be a free multiplicative
abelian group generated by $x_{n+1},\cdots,x_{n+m}$. We set ${\mathbb
Z}\cP:={\mathbb Z}[x_{n+1}^{\pm1}, \cdots ,x_{n+m}^{\pm1}]$. Let
$K:=\{\frac{g}{h} |\ g,\ h \in {\mathbb Z}\cP,\ h\neq0 \}$ be the field
of fractions of ${\mathbb Z}\cP$, and $\cF:=K(x_{1}, \cdots ,x_{n})$ 
be the field of rational functions.

\subsection{Cluster algebras of geometric type}
\begin{defn}
Let $B=(b_{ij})$ be an $n\times n$ integer matrix. 
\begin{enumerate}
\item $B$ is {\it skew symmetric} if $b_{ij}=-b_{ji}$ for any $i,\ j \in [1,n]$.
\item $B$ is  {\it skew symmetrizable} if there exists a positive 
integer diagonal matrix $D$ such that $DB$ is skew symmetric.
\item $B$ is  {\it sign skew symmetric} if $b_{ij}b_{ji}\leq0$ for any
$i,\ j \in [1,n]$, 
and if $b_{ij}b_{ji}=0$ then $b_{ij}=b_{ji}=0$.
\end{enumerate}
\end{defn}

Note that each skew symmetric matrix is skew symmetrizable, and each skew symmetrizable matrix is sign skew symmetric.

\begin{defn}
We set $n$-tuple of variables $\textbf{x}=(x_1,\cdots,x_n)$. Let $\tilde{B}=(b_{ij})_{1\leq i\leq n+m,\ 1\leq j \leq n}$ be $(n+m)\times
n$ integer matrix whose principal part $B:=(b_{ij})_{1\leq i,j\leq n}$ 
is sign skew symmetric. Then a pair $\Sigma=(\textbf{x},\tilde{B})$ is called a {\it seed}, \textbf{x} a cluster and $x_1, \cdots ,\ x_n$ cluster variables. For a seed $\Sigma=(\textbf{x},\tilde{B})$, principal part $B$ of $\tilde{B}$ is called the {\it exchange matrix}. 
\end{defn}
\begin{defn}
If $B$ is skew symmetric (resp. skew symmetrizable, sign skew symmetric),
 we say $\tilde{B}$ is skew symmetric
(resp. skew symmetrizable, sign skew symmetric).
\end{defn}
\begin{defn}\label{adc}
For a seed $\Sigma=(\textbf{x}, \tilde{B}=(b_{ij}))$,  
an {\it adjacent cluster} in direction $k\in [1,n]$ is defined by 
\[ \textbf{x}_k = (\textbf{x}\setminus \{x_k\})\cup \{x_k'\} ,\]
where $x_k'$ is the new cluster variable defined by the  {\it exchange relation}
\[ x_k x_k' = 
\prod_{1\leq i \leq n+m,\ b_{ik}>0} x_i^{b_{ik}}
+\prod_{1\leq i \leq n+m,\ b_{ik}<0} x_i^{-b_{ik}}. \]
\end{defn}

\begin{defn}
Let $A=(a_{ij}),\ A'=(a_{ij}')$ be two matrices of the same size. We say
 that $A'$ is obtained from $A$ by the matrix mutation in direction $k$,
 and denote $A'=\mu_k(A)$ if
\[
a_{ij}'=\begin{cases}
		-a_{ij} & {\rm if}\ i=k\ {\rm or}\ j=k, \\
		a_{ij}+\frac{|a_{ik}|a_{kj}+a_{ik}|a_{kj}|}{2} & {\rm otherwise}.
		\end{cases}
\]
For $A,\ A'$, if there exists a finite sequence $(k_1,\cdots,k_s)$,
$(k_i\in[1,n])$ such that $A'=\mu_{k_1}\cdots\mu_{k_s}(A)$, we say $A$
is mutation equivalent to $A'$, 
and denote $A \cong A'$.
\end{defn}

\begin{prop}${\cite{M-M-A}}$
For $k\in [1,n]$, $\mu_k(\mu_k(A))=A$.
\end{prop}

\begin{defn}
Let $A$ be a sign skew symmetric matrix. We say $A$ is {\it totally sign skew symmetric} if any matrix that is mutation equivalent to $A$ is sign skew symmetric. Then a seed $(\textbf{x},A)$ is called a totally mutable seed.
\end{defn}

Next proposition can be easily verified by the definition of $\mu_k$:

\begin{prop}\label{totally}${\cite[Proposition\ 3.6]{M-M-A}}$
Skew symmetrizable matrices are totally sign skew symmetric.
\end{prop}

For a seed $\Sigma=(\textbf{x},\tilde{B})$, we say that the seed $\Sigma'=(\textbf{x}',\tilde{B}')$ is adjacent to $\Sigma$ if $\textbf{x}'$ is adjacent to $\textbf{x}$ in direction $k$ and $\tilde{B}'=\mu_k(\tilde{B})$. Two seeds $\Sigma$ and $\Sigma_0$ are mutation equivalent if one of them can be obtained from another seed by a sequence of pairwise adjacent seeds and we denote $\Sigma \sim\Sigma_0$.

Now let us define the cluster algebra of geometric type.

\begin{defn}
Let $\tilde{B}$ be a skew symmetrizable matrix, and $\Sigma=(\textbf{x},\tilde{B})$ a seed. The cluster algebra (of geometric type) $\cA=\cA(\Sigma)$ associated with
seed $\Sigma$ is defined as the ${\mathbb Z}\cP$-subalgebra of $\cF$ generated by all cluster variables in all seeds which are mutation equivalent to $\Sigma$.
\end{defn}

For a seed $\Sigma$, we define ${\mathbb Z}\cP$-subalgebra $\UU(\Sigma)$ of $\cF$ by
\[ \UU(\Sigma):={\mathbb Z}\cP[\textbf{x}^{\pm 1}] \cap {\mathbb Z}\cP[\textbf{x}_1^{\pm 1}] \cap \cdots \cap {\mathbb Z}\cP[\textbf{x}_n^{\pm 1}] .\]
 Here, ${\mathbb Z}\cP[\textbf{x}^{\pm 1}]$ is the Laurent polynomial ring in \textbf{x}.

\begin{defn}\label{upper}
Let $\Sigma_0$ be a totally mutable seed. We define an {\it upper cluster algebra} $\ovl{\cA}=\ovl{\cA}(\Sigma_0)$ as the intersection of the subalgebras $\UU(\Sigma)$ for all seeds $\Sigma \sim\Sigma_0$.
\end{defn}

For a totally mutable seed $\Sigma$, following the inclusion relation holds \cite{A-F-Z}:
\[ \cA(\Sigma) \subset \ovl{\cA}(\Sigma) .\]

\subsection{Upper cluster algebra $\ovl{\cA}(\textbf{i})$}

As in Sect.\ref{DBCs}, let $G=Sp_{2r}(\mathbb{C})$ be the simple algebraic group of type ${\rm C}_r$ and $W$ be its Weyl group. We set $u\in W$ and its reduced word $\textbf{i}$ as in (\ref{uset0}) and (\ref{iset0}):
\begin{equation}\label{uset00}
u=\underbrace{s_1s_2\cdots s_r}_{1{\rm \,st\
cycle}}
\underbrace{s_1\cdots s_{r}}_{2{\rm \,nd\
cycle}}
\cdots \underbrace{s_1\cdots s_{r}}_{m-1{\rm \,\,th\ cycle}}
\underbrace{s_1\cdots s_{i_n}}_{m{\rm \,\,th\ cycle}}, 
\end{equation}
\begin{equation}\label{iset00}
\textbf{i}=(\underbrace{1,\cdots,r}_{1{\rm \,st\
cycle}},\underbrace{1,\cdots,r}_{2{\rm \,nd\
cycle}},\cdots,\underbrace{1,\cdots,r}_{m-1{\rm \,\,th\ cycle}},
\underbrace{1,\cdots,i_n}_{m{\rm \,\,th\ cycle}}). 
\end{equation}

In this subsection, we define the upper cluster algebra $\ovl{\cA}(\textbf{i})$, which satisfies that $\ovl{\cA}(\textbf{i})\otimes \mathbb{C}$ is isomorphic to the coordinate ring $\mathbb{C}[G^{u,e}]$ of the double Bruhat cell \cite{A-F-Z}. Let $i_k$ $(k\in[1,l(u)])$ be the $k$-th index of $\textbf{i}$ from the left. 

At first, we define a set e(\textbf{i}) as 
\[ e(\textbf{i}):= \{k| {\rm There\ exist\ some}\ l>k\ {\rm such\ that}\ i_k=i_l \}.  \]

Next, let us define a matrix $\tilde{B}=\tilde{B}(\textbf{i})$. 

\begin{defn}
Let $\tilde{B}(\textbf{i})$ be an integer matrix with rows labelled by all the indices in $[-1,-r]\cup [1,l(u)]$ and columns labelled by all the indices in $e(\textbf{i})$. For $k\in[-1,-r]\cup [1,l(u)]$ and $l\in e(\textbf{i})$, an entry $b_{kl}$ of $\tilde{B}(\textbf{i})$ is determined as follows: 
\[
b_{kl}=\begin{cases}
		-{\rm sgn}\left((k-l)\cdot i_p\right) & {\rm if}\ p=q, \\
		-{\rm sgn}((k-l)\cdot i_p\cdot a_{|i_k||i_l|}) & {\rm if}\ p<q\ {\rm and}\ {\rm sgn}(i_p\cdot i_q)(k-l)(k^+-l^+)>0,  \\
		0 & {\rm otherwise}.
	\end{cases}
\]
\end{defn}

\begin{prop}\label{propss}${\cite[Proposition\ 2.6]{A-F-Z}}$
$\tilde{B}(\textbf{i})$ is skew symmetrizable. 
\end{prop}

By Proposition $\ref{totally}$, Definition $\ref{upper}$ and Proposition \ref{propss}, we can construct the upper cluster algebra:

\begin{defn}
We denote this upper cluster algebra by $\ovl{\mathcal{A}}(\textbf{i})$.
\end{defn}

\subsection{Generalized minors and bilinear form}\label{gmtypeA}

As in the previous section, we set $G=Sp_{2r}(\mathbb{C})$, $u\in W$ and its reduced word $\textbf{i}$ as in (\ref{uset00}) and (\ref{iset00}). We also set $\bar{\cA}(\textbf{i})_{\mathbb{C}}:=\bar{\cA}(\textbf{i})\otimes \mathbb{C}$ and  $\cF_{\mathbb{C}}:=\cF \otimes \mathbb{C}$. It is known that the coordinate ring $\mathbb{C}[G^{u,e}]$ of the double Bruhat cell is isomorphic to $\bar{\cA}(\textbf{i})_{\mathbb{C}}$ (Theorem \ref{clmainthm}). To describe this isomorphism explicitly, we need generalized minors.  

We set $G_0:=N_-HN$, and let $x=[x]_-[x]_0[x]_+$ with $[x]_-\in N_-$, $[x]_0\in H$, $[x]_+\in N$ be the corresponding decomposition. 

\begin{defn}
For $i\in[1,r]$ and $w$, $w'\in W$, the generalized minor $\Delta_{w\Lambda_i,w'\Lambda_i}$ is a regular function on $G$ whose restriction to the open set $wG_0w'^{-1}$ is given by $\Delta_{w\Lambda_i,w'\Lambda_i}(x)=([w^{-1}x w']_0)^{\Lambda_i}$. Here, $\Lambda_i$ is the $i$-th  fundamental weight. In particular, we write $\Delta_{\Lambda_i}:=\Delta_{\Lambda_i,\Lambda_i}$ and call it {\it principal minor}.
\end{defn}

We set $\ge={\rm Lie}(G)$. Let $\omega:\ge\to\ge$ be the anti involution 
\[
\omega(e_i)=f_i,\q
\omega(f_i)=e_i,\q \omega(h)=h,
\] and extend it to $G$ by setting
$\omega(x_i(c))=y_{i}(c)$, $\omega(y_{i}(c))=x_i(c)$ and $\omega(t)=t$
$(t\in H)$. Here, $x_i$ and $y_i$ were defined in Sect.\ref{factproC} (\ref{xiyidef}).

There exists a $\ge$ (or $G$)-invariant bilinear form on the
finite-dimensional  irreducible
$\ge$-module $V(\lm)$ such that 
\[
 \lan au,v\ran=\lan u,\omega(a)v\ran,
\q\q(u,v\in V(\lm),\,\, a\in \ge\ (\text{or }G)).
\]
For $g\in G$, 
we have the following simple fact:
\[
 \Del_{\Lm_i}(g)=\lan gu_{\Lm_i},u_{\Lm_i}\ran,
\]
where $u_{\Lm_i}$ is a properly normalized highest weight vector in
$V(\Lm_i)$. Hence, for $w,w'\in W$, we have
\begin{equation}\label{minor-bilin}
 \Del_{w\Lm_i,w'\Lm_i}(g)=
\Del_{\Lm_i}({\ovl w}^{-1}g\ovl w')=
\lan {\ovl w}^{-1}g\ovl w'\cdot u_{\Lm_i},u_{\Lm_i}\ran
=\lan g\ovl w'\cdot u_{\Lm_i}\, ,\, \ovl{w}\cdot u_{\Lm_i}\ran,
\end{equation}
where $\ovl w$ is the one we defined in Sect.\ref{factpro} (\ref{smpl}), and note that $\omega(\ovl s_i^{\pm})=\ovl s_i^{\mp}$.

\subsection{Cluster algebras on Double Bruhat cells of type C}\label{cdb}

For $k \in [1,l(u)]$, let  $i_k$ be the $k$-th index of $\textbf{i}$ $(\ref{iset00})$ from the left, and we suppose that it belongs to the $m'$ th cycle. We set
\begin{equation}\label{inc}
u_{\leq k}=u_{\leq k}(\textbf{i}):=\underbrace{s_1s_2\cdots s_r}_{1{\rm \,st\
cycle}}
\underbrace{s_1\cdots s_{r}}_{2{\rm \,nd\
cycle}}
\cdots \underbrace{s_1\cdots s_{i_k}}_{m'{\rm \,th\
cycle}}.
\end{equation}
For $k \in [-1,-r]$, we set $u_{\leq k}:=e$ and $i_k:=k$. For $k \in [-1,-r]\cup [1,\ l(u)]$, we define
\[ \Delta(k;\textbf{i})(x):=\Delta_{u_{\leq k} \Lambda_{i_k},\Lambda_{i_k}}(x).  \]

Finally, we set
\[ F(\textbf{i}):=\{ \Delta(k;\textbf{i})(x)|k \in [-1,-r]\cup[1,\ l(u)] \}. \]
It is known that the set $F(\textbf{i})$ is an algebraically independent generating set for the field of rational functions $\mathbb{C}(G^{u,e})$ \cite[Theorem 1.12]{F-Z}. Then, we have the following theorem.

\begin{thm}\label{clmainthm}${\cite[Theorem\ 2.10]{A-F-Z}}$
The isomorphism of fields $\varphi :F_{\mathbb{C}} \rightarrow \mathbb{C}(G^{u,e})$ defined by $\varphi (x_k)=\Delta(k;\textbf{i})\ (k \in [-1,-r]\cup [1,l(u)] )$ restricts to an isomorphism of algebras $\bar{\cA}(\textbf{i})_{\mathbb{C}}\rightarrow \mathbb{C}[G^{u,e}]$.
\end{thm}

\section{Explicit formulas of cluster variables}\label{gmc}

In the rest of the paper, we consider the case
$G=Sp_{2r}(\mathbb{C})$. 
Let $u\in W$ be
\begin{equation}\label{uvset}
u=(s_1s_2\cdots s_r)^{m-1}s_1\cdots s_{i_n},
\end{equation} 
where $n=l(u)$, $1\leq i_n\leq r$ and $1\leq m\leq r$. Let
\begin{equation}\label{iset}
 \textbf{i}=(\underbrace{1,\cdots,r}_{1{\rm \,st\
  cycle}},\underbrace{1,\cdots,r}_{2{\rm \,nd\
  cycle}},\cdots,\underbrace{1,\cdots,r}_{m-1{\rm \,th\ cycle}},
\underbrace{1,\cdots,i_n}_{m{\rm \,th\ cycle}}), 
\end{equation}
be a reduced word $\textbf{i}$ for $u$, that is, \textbf{i} is the left factor of $(1,2,3,\cdots,r)^r$.  Let $i_k$ be the $k$-th index of  $\textbf{i}$ from the left, and belong to $m'$-th cycle. As we shall show in lemma \ref{gmlem}, we may assume $i_n=i_k$. 

By Theorem \ref{clmainthm}, we can regard $\mathbb{C}[G^{u,e}]$ as an
upper cluster algebra and $\{\Delta(k;\textbf{i})\}$ as its cluster
variables. Each $\Delta(k;\textbf{i})$ is a regular function on
$G^{u,e}$. On the other hand, by Proposition \ref{gprime} (resp. Theorem
\ref{fp2}), we can consider $\Delta(k;\textbf{i})$ as a function on
$H\times
(\mathbb{C}^{\times})^{l(u)}$ (resp. $(\mathbb{C}^{\times})^{l(u)}$). Then
we change the variables of 
$\{\Delta(k;\textbf{i})\}$ as follows: 
\begin{defn}\label{gendef}For $a\in H$ and
\begin{multline}\label{yset}
 \textbf{Y}:=(Y_{1,1},Y_{1,2},\cdots,Y_{1,r},Y_{2,1},Y_{2,2},\cdots,Y_{2,r}, \\
\cdots,Y_{m-1,1},\cdots,Y_{m-1,r},Y_{m,1},\cdots,Y_{m,i_n})\in  (\mathbb{C}^{\times})^{n}, 
\end{multline}
we set 
\[ 
\Delta^G(k;\textbf{i})(a,\textbf{Y}):=(\Delta(k;\textbf{i})\circ
 \ovl{x}^G_{\textbf{i}})(a,\textbf{Y}), 
\] 
\[ 
\Delta^L(k;\textbf{i})(\textbf{Y}):=(\Delta(k;\textbf{i})\circ
 x^L_{\textbf{i}})(\textbf{Y}).
\]
\end{defn}
We will describe the function $\Delta^L(k;\textbf{i})(\textbf{Y})$ explicitly since $\Delta^G(k;\textbf{i})(a,\textbf{Y})$ is immediately obtained from $\Delta^L(k;\textbf{i})(\textbf{Y})$ (Proposition \ref{gprop}).

\begin{rem}\label{importantrem}
If we see the variables $Y_{s,0}$, $Y_{s,r+1}$ 
$(1\leq s\leq m)$ then 
we understand
\[ Y_{s,0}=Y_{s,r+1}=1. \]
For example, if $i=1$ then
\[ Y_{s,i-1}=1. \]
\end{rem}

\subsection{Generalized minor $\Delta^G(k;\textbf{i})(a,\textbf{Y})$}

In this subsection, we shall prove that
$\Delta^G(k;\textbf{i})(a,\textbf{Y})$ is 
immediately obtained from $\Delta^L(k;\textbf{i})(\textbf{Y})$:
\begin{prop}\label{gprop}
We set $d:=i_k$. For $a=t^{\sum_{i}a_i h_i}\in H$ $(t\in \mathbb{C}^{\times})$, we have
\[ 
\Delta^G(k;\textbf{i})(a,\textbf{Y})=
\begin{cases}
t^{(a_{r}-a_{m'}-a_{d-r+m'})}\Delta^L(k;\textbf{i})(\textbf{Y}) & {\rm if}\ m'+d>r, \\
t^{(a_{m'+d}-a_{m'})}\Delta^L(k;\textbf{i})(\textbf{Y}) & {\rm if}\ m'+d\leq r.
\end{cases}
\] 
\end{prop}
This proposition follows from (\ref{wtv}) and the following lemma:
\begin{lem}\label{gllem}
In the above setting, if $m'+d>r$, then we have
\[ \Delta^G(k;\textbf{i})(a,\textbf{Y})=\lan  ax^L_{\textbf{i}}(\textbf{Y})(v_1\wedge v_2\wedge\cdots \wedge v_d) ,\q v_{m'+1}\wedge\cdots\wedge v_r\wedge v_{\ovl{d-r+m'}}\wedge\cdots\wedge v_{\ovl{1}}  \ran, \]
\[ \Delta^L(k;\textbf{i})(\textbf{Y})=\lan  x^L_{\textbf{i}}(\textbf{Y})(v_1\wedge v_2\wedge\cdots \wedge v_d) ,\q v_{m'+1}\wedge\cdots\wedge v_r\wedge v_{\ovl{d-r+m'}}\wedge\cdots\wedge v_{\ovl{1}}  \ran, \]
where $\lan ,\ran$ is the bilinear form we defined in Sect.\ref{gmtypeA}. In the case $m'+d\leq r$, we have
\[ \Delta^G(k;\textbf{i})(a,\textbf{Y})=\lan  ax^L_{\textbf{i}}(\textbf{Y})(v_1\wedge v_2\wedge\cdots \wedge v_d) ,\q v_{m'+1}\wedge\cdots\wedge v_{m'+d}  \ran, \]
\[ \Delta^L(k;\textbf{i})(\textbf{Y})=\lan  x^L_{\textbf{i}}(\textbf{Y})(v_1\wedge v_2\wedge\cdots \wedge v_d) ,\q v_{m'+1}\wedge\cdots\wedge v_{m'+d}  \ran. \]
\end{lem}
\nd
{\sl Proof.}
Let us prove this lemma for $\Delta^L(k;\textbf{i})(\textbf{Y})$ since the
case for $\Delta^G(k;\textbf{i})(a,\textbf{Y})$ 
is proven similarly.
Using (\ref{minor-bilin}) and (\ref{inc}), we see that $\Delta^L(k;\textbf{i})(\textbf{Y})=\Delta_{u_{\leq k}\Lambda_d, \Lambda_d}(x^L_{\textbf{i}}(\textbf{Y}))$ is given as 
\begin{equation}\label{bilin2}
\lan x^L_{\textbf{i}}(\textbf{Y})(v_1\wedge v_2\wedge\cdots \wedge v_d) ,\q \underbrace{\ovl{s_1}\cdots \ovl{s_r}}_{1{\rm \,st\
 cycle}}\cdots\underbrace{\ovl{s_1}\cdots \ovl{s_d}}_{m'{\rm \,th\ cycle}} (v_1\wedge v_2\wedge\cdots \wedge v_d) \ran.
\end{equation}
By (\ref{smpl}), for $1\leq i \leq r-1$ and $1\leq j \leq r$, we get
\[ \ovl{s_i} v_j=
\begin{cases}
v_{i+1} & {\rm if}\ j=i, \\
-v_{i} & {\rm if}\ j=i+1, \\
v_j & {\rm if}\ {\rm otherwise},
\end{cases}\q 
\ovl{s_i} v_{\ovl{j}}=
\begin{cases}
v_{\ovl{i}} & {\rm if}\ j=i+1, \\
-v_{\ovl{i+1}} & {\rm if}\ j=i, \\
v_{\ovl{j}} & {\rm if}\ {\rm otherwise},
\end{cases}
\]
and we obtain
\[ \ovl{s_r} v_j=
\begin{cases}
v_{\ovl{r}} & {\rm if}\ j=r, \\
v_j & {\rm if}\ j\neq r,
\end{cases}\q 
\ovl{s_r} v_{\ovl{j}}=
\begin{cases}
-v_{r} & {\rm if}\ j=r, \\
v_{\ovl{j}} & {\rm if}\ j\neq r.
\end{cases}
\]
Therefore, if $m'+d\leq r$, then
\begin{equation}\label{ukeq1}u_{\leq k}(v_1\wedge\cdots\wedge v_d)=
\underbrace{\ovl{s_1}\cdots \ovl{s_r}}_{1{\rm \,st\ 
 cycle}}\cdots\underbrace{\ovl{s_1}\cdots \ovl{s_d}}_{m'{\rm \,th\ cycle}} (v_1\wedge \cdots \wedge v_d)=v_{m'+1}\wedge v_{m'+2}\wedge\cdots \wedge v_{m'+d}. 
\end{equation}
If $m'+d> r$, then we get
\begin{eqnarray}\label{ukeq2}
& &u_{\leq k}(v_1\wedge\cdots\wedge v_d)\nonumber \\
&=&
\underbrace{\ovl{s_1}\cdots \ovl{s_r}}_{1{\rm \,st\
 cycle}}\cdots\underbrace{\ovl{s_1}\cdots \ovl{s_d}}_{m'{\rm \,th\ cycle}} (v_1\wedge v_2\wedge\cdots \wedge v_d) \nonumber \\
&=&\underbrace{\ovl{s_1}\cdots \ovl{s_r}}_{1{\rm \,st
\ cycle}}\cdots\underbrace{\ovl{s_1}\cdots \ovl{s_r}}_{m'-r+d{\rm \,th\ cycle}} (v_{r-d+1}\wedge\cdots \wedge v_{r}) \nonumber \\
&=&\underbrace{\ovl{s_1}\cdots \ovl{s_r}}_{1{\rm \,st
\ cycle}}\cdots\underbrace{\ovl{s_1}\cdots \ovl{s_r}}_{m'-r+d-1{\rm \,th\ cycle}} (v_{r-d+2}\wedge\cdots \wedge v_{r}\wedge v_{\ovl{1}})\nonumber \\
&=&\underbrace{\ovl{s_1}\cdots \ovl{s_r}}_{1{\rm \,st
\ cycle}}\cdots\underbrace{\ovl{s_1}\cdots \ovl{s_r}}_{m'-r+d-2{\rm \,th\ cycle}} (v_{r-d+3}\wedge\cdots \wedge v_{r}\wedge v_{\ovl{1}}\wedge (-v_{\ovl{2}}))\nonumber \\
&=&\cdots\q =v_{m'+1}\wedge\cdots \wedge v_{r}\wedge v_{\ovl{1}}\wedge (-v_{\ovl{2}})\wedge 
 ((-1)^2v_{\ovl{3}})\wedge\cdots\wedge ((-1)^{d-r+m'-1}v_{\ovl{d-r+m'}})
 \nonumber \\
&=&v_{m'+1}\wedge\cdots \wedge v_{r}\wedge v_{\ovl{d-r+m'}}\wedge\cdots\wedge v_{\ovl{1}}.
\end{eqnarray}
Hence, we get our claim by (\ref{bilin2}). \qed 

In the rest of the paper, we will treat $\Delta^L(k;\textbf{i})(\textbf{Y})$ only by Proposition \ref{gprop}.

\subsection{Generalized minor $\Delta^L(k;\textbf{i})(\textbf{Y})$}\label{Mainsec}

\begin{lem}\label{gmlem}
Let $u$, $\textbf{i}$ and $\textbf{Y}$ be as in $(\ref{uvset})$, $(\ref{iset})$ and $(\ref{yset})$. Let $i_{n+1}\in [1,r]$ be an index such that $u':=us_{i_{n+1}}\in W$ satisfies $l(u')>l(u)$. We set the reduced word $\textbf{i}'$ for $u'$ as
\[
\textbf{i}'=(\underbrace{1,\cdots,r}_{1{\rm \,st\
  cycle}},\underbrace{1,\cdots,r}_{2{\rm \,nd\
  cycle}},\cdots,\underbrace{1,\cdots,r}_{m-1{\rm \,th\ cycle}},
\underbrace{1,\cdots,i_n}_{m{\rm \,th\ cycle}},i_{n+1}),
\]
and denote $\textbf{Y}'\in(\mathbb{C}^{\times})^{n+1}$ by
\[ \textbf{Y}':=(Y_{1,1},\cdots,Y_{1,r},\cdots,Y_{m-1,1},\cdots,Y_{m-1,r},Y_{m,1},\cdots,Y_{m,i_n},Y).
\] 
For an integer $k$ $(1\leq k\leq n)$, if $d:=i_k\neq i_{n+1}$, then $\Delta^L(k;\textbf{i}')(\textbf{Y}')$ does not depend on $Y$, so we can regard it as a function on $(\mathbb{C}^{\times})^{n}$. Furthermore, we have 
\begin{equation}\label{Lomit}
\Delta^L(k;\textbf{i})(\textbf{Y})=\Delta^L(k;\textbf{i}')(\textbf{Y}') .
\end{equation}
\end{lem}
\nd
{\sl Proof.}
By the definition (\ref{xldef}) of $x^L_{\textbf{i}}$, we have
\begin{equation}\label{gmlempr1}
x^L_{\textbf{i}'}(\textbf{Y})=x^L_{\textbf{i}}(\textbf{Y})x_{-i_{n+1}}(Y). 
\end{equation}
On the other hand, since $f_{i}^2=0$ on $V(\Lambda_1)$, we have ${\rm exp}(t f_{i})=1+t f_i$ $(i=1,\cdots,r,\ t\in \mathbb{C})$. Hence, by $x_{-i_{n+1}}(Y):={\rm exp}(Y f_{i_{n+1}})\cdot(Y^{-h_{i_{n+1}}})$ (see (\ref{alxmdef})), we get
\begin{equation}\label{xmpro}
x_{-i_{n+1}}(Y)v_j=
\begin{cases}
Y^{-1}v_{i_{n+1}}+v_{i_{n+1}+1} & {\rm if}\ j=i_{n+1}, \\
Y v_{i_{n+1}+1} & {\rm if}\ j=i_{n+1}+1, \\
v_j & {\rm otherwise},
\end{cases}
\end{equation}
where in the case $j=i_{n+1}$, we set $v_{r+1}:=v_{\ovl{r}}$. Thus, if $d<i_{n+1}$, then we have $x_{-i_{n+1}}(Y)(v_1\wedge\cdots\wedge v_d)=v_1\wedge\cdots\wedge v_d$.
If $d>i_{n+1}$, then we have
\begin{flushleft}
$x_{-i_{n+1}}(Y)(v_1\wedge\cdots\wedge v_d)$
\end{flushleft}
\begin{eqnarray*}
&=&v_1\wedge\cdots \wedge v_{i_{n+1}-1}\wedge (Y^{-1}v_{i_{n+1}}+v_{i_{n+1}+1})\wedge Y v_{i_{n+1}+1} \wedge\cdots \wedge v_d \\
&=&v_1\wedge\cdots\wedge v_d.
\end{eqnarray*}
Since we assume $i_{n+1}\neq d$, we get
\begin{equation}\label{gmlempr2}
x_{-i_{n+1}}(Y)(v_1\wedge\cdots\wedge v_d)=v_1\wedge\cdots\wedge v_d.
\end{equation}
We can easily see that $u_{\leq k}=u'_{\leq k}(=\underbrace{s_1\cdots s_r}_{1{\rm \,st\
 cycle}}\cdots\underbrace{s_1\cdots s_d}_{m'{\rm \,th\ cycle}})$. Therefore, it follows from (\ref{minor-bilin}), (\ref{gmlempr1}) and (\ref{gmlempr2}) that
\begin{eqnarray*}
\Delta^L(k;\textbf{i}')(\textbf{Y}')&=&\Delta_{u'_{\leq k}\Lambda_d,\Lambda_d}(x^L_{\textbf{i}'}(\textbf{Y}')) \\
&=&\lan x^L_{\textbf{i}'}(\textbf{Y}')(v_1\wedge v_2\wedge\cdots \wedge v_d) ,\ u'_{\leq k} (v_1\wedge v_2\wedge\cdots \wedge v_d) \ran \\
&=&\lan x^L_{\textbf{i}}(\textbf{Y})(v_1\wedge v_2\wedge\cdots \wedge v_d) ,\ u_{\leq k} (v_1\wedge v_2\wedge\cdots \wedge v_d) \ran=\Delta^L(k;\textbf{i})(\textbf{Y}),
\end{eqnarray*}
which is our desired result. \qed

\vspace{5mm}

By this lemma, when we calculate $\Delta^L(k;\textbf{i})(\textbf{Y})$, we may
assume that $i_n=i_k$ 
without loss of generality.

For $1\leq l\leq m$ and $1\leq k\leq r$, we set the Laurent monomials
\begin{equation}\label{ccbar}
 \ovl{C}(l,k):=\frac{Y_{l,k-1}}{Y_{l,k}} ,\q C(l,k):=\frac{Y_{l,k+1}}{Y_{l+1,k}}.
\end{equation}

\begin{rem}
In \cite{N1}, it was defined $\ovl{C}^{(l)}_k:=\frac{Y_{r-l,k-1}}{Y_{r-l,k}}$ and $C^{(l)}_k:=\frac{Y_{r-l,k}}{Y_{r-l+1,k-1}}$, which coincide with $\ovl{C}(r-l,k)$ and $C(r-l,k-1)$ in $(\ref{ccbar})$ respectively.   
\end{rem}

For $1\leq l\leq r$, we set $|l|=|\ovl{l}|=l$. The following theorem is our main result. 

\begin{thm}\label{thm1}
In the above setting, we set $d:=i_k=i_n$ and
\[ \textbf{Y}:=(Y_{1,1},Y_{1,2},\cdots,Y_{1,r},\cdots,Y_{m-1,1},\cdots,Y_{m-1,r},Y_{m,1},\cdots,Y_{m,i_n})\in  (\mathbb{C}^{\times})^{n}. \]
Then we have
\begin{multline*} 
\Delta^L(k;\textbf{i})(\textbf{Y})= \sum_{(*)}\prod^{d}_{i=1} \ovl{C}(m-l^{(1)}_i,k^{(1)}_i)\cdots \ovl{C}(m-l^{(\delta_i)}_i,k^{(\delta_i)}_i)\\
\cdot C(m-l^{(\delta_i+1)}_i,|k^{(\delta_i+1)}_i|-1)\cdots C(m-l^{(m-m')}_i,|k^{(m-m')}_i|-1),
\end{multline*} 
\[ l^{(s)}_{i}:=
\begin{cases}
k^{(s)}_{i}+s-i-1 & {\rm if}\  s\leq \del_{i}, \\
s-i+r & {\rm if}\  s>\del_{i}
\end{cases}
\q (1\leq i\leq d)
\]
where $(*)$ is the conditions for $k^{(s)}_i$ $(1\leq s\leq m-m',\ 1\leq i\leq d)$ : 
$1\leq k^{(s)}_1<k^{(s)}_2<\cdots<k^{(s)}_d\leq\ovl{1}\q (1\leq s\leq m-m')$,
$1\leq k^{(1)}_i\leq \cdots\leq k^{(m-m')}_i\leq m'+i\q (1\leq i\leq r-m')$,
and $1\leq k^{(1)}_i\leq \cdots\leq k^{(m-m')}_i\leq \ovl{1}\q (r-m'+1\leq i\leq d)$, and 
$\delta_i$ $(i=1,\cdots,d)$ are the numbers which satisfy $1\leq k^{(1)}_i\leq k^{(2)}_i\leq \cdots\leq k^{(\delta_i)}_i\leq r$, $\ovl{r}\leq k^{(\delta_i+1)}_i\leq \cdots\leq k^{(m-m')}_i\leq\ovl{1}$.
\end{thm}

\begin{ex}\label{pathex3}
For {\rm rank} $r=3$, $u=s_1s_2s_3s_1s_2s_3s_1s_2$, $k=5$ 
and the reduced word $\textbf{i}=(-1,-2,-3,-1,-2,-3,-1,-2)$ for $u$, 
we have $m=3$, $m'=2$ and $d=2$ $($see $(\ref{uvset}),\ (\ref{iset}))$. 
Then, we have $s=1$ and write $k_i$ for $k^{(s)}_i$. 
Thus, the set of  all $(k_1,k_2)$ satisfying $(*)$ in Theorem \ref{thm1} is
\[
\{(1,2),(1,3),(1,\ovl 3),(1,\ovl 2),(1,\ovl 1),(2,3),(2,\ovl 3),(2,\ovl 2),
(2,\ovl 1),(3,\ovl 3),(3,\ovl 2),(3,\ovl 1)\}
\]
Here, for all $(k_1,k_2)$ the corresponding monomials are as follows:
\[
\begin{array}{ccc}
(1,2)\leftrightarrow  \ovl{C}(3,1)\ovl{C}(3,2)&
(1,3)\leftrightarrow    \ovl{C}(3,1)\ovl{C}(2,3)&
(1,\ovl 3)\leftrightarrow \ovl{C}(3,1)C(1,2)\\
(1,\ovl 2)\leftrightarrow \ovl{C}(3,1)C(1,1) &
(1,\ovl 1)\leftrightarrow \ovl{C}(3,1)C(1,0)&
(2,3)\leftrightarrow   \ovl{C}(2,2)\ovl{C}(2,3)\\
(2,\ovl 3)\leftrightarrow  \ovl{C}(2,2)C(1,2)&
(2,\ovl 2)\leftrightarrow \ovl{C}(2,2)C(1,1) &
(2,\ovl 1)\leftrightarrow \ovl{C}(2,2)C(1,0)\\
(3,\ovl 3)\leftrightarrow \ovl{C}(1,3)C(1,2)&
(3,\ovl 2)\leftrightarrow  \ovl{C}(1,3)C(1,1)&
(3,\ovl 1)\leftrightarrow \ovl{C}(1,3)C(1,0)
\end{array}
\]
Thus, we obtain:
\begin{eqnarray*}
\Delta^L(5;\textbf{i})(\textbf{Y})&=&\ovl{C}(3,1)\ovl{C}(3,2)
+\ovl{C}(3,1)\ovl{C}(2,3)
+\ovl{C}(3,1)C(1,2)
+\ovl{C}(3,1)C(1,1) \\
& &
+\ovl{C}(3,1)C(1,0)
+\ovl{C}(2,2)\ovl{C}(2,3)
+\ovl{C}(2,2)C(1,2)
+\ovl{C}(2,2)C(1,1)\\
& &+\ovl{C}(2,2)C(1,0)
+\ovl{C}(1,3)C(1,2)
+\ovl{C}(1,3)C(1,1)
+\ovl{C}(1,3)C(1,0)   \\
&=&
\frac{1}{Y_{3,2}}+\frac{Y_{2,2}}{Y_{3,1}Y_{2,3}}+\frac{Y_{1,3}}{Y_{3,1}Y_{2,2}}
+\frac{Y_{1,2}}{Y_{3,1}Y_{2,1}} 
+\frac{Y_{1,1}}{Y_{3,1}}+\frac{Y_{2,1}}{Y_{2,3}}+\frac{Y_{2,1}Y_{1,3}}{Y_{2,2}^2}\\
&&+2\frac{Y_{1,2}}{Y_{2,2}} 
+\frac{Y_{2,1}Y_{1,1}}{Y_{2,2}}
+\frac{Y_{1,2}^2}{Y_{2,1}Y_{1,3}}+\frac{Y_{1,1}Y_{1,2}}{Y_{1,3}}.
\end{eqnarray*}
Note that since $\ovl{C}(1,3)C(1,2)=\ovl{C}(2,2)C(1,1)=\frac{Y_{1,2}}{Y_{2,2}}$, 
the coefficient of $\frac{Y_{1,2}}{Y_{2,2}}$ in the above formula 
is equal to $2$.
\end{ex}

\section{The proof of Theorem \ref{thm1}}\label{prsec}

In this section, we shall give the proof of Theorem \ref{thm1}.

\subsection{The set $X_d(m,m')$ of paths: path descriptions}

In this subsection, we shall introduce a set $X_d(m,m')$ of ``paths''
which correspond to the terms of $\Delta^L(k;\textbf{i})(\textbf{Y})$, 
which we call {\it path descriptions} of generalized minors. Let
$m$, $m'$ and $d$ be the positive integers as in 
\ref{Mainsec}. We set $J:=\{j,\ovl{j}|\ 1\leq j\leq r\}$ and for $1\leq l\leq r$, set $|l|=|\ovl{l}|=l$.
\begin{defn}
Let us define the directed graph $(V_d,E_d)$ as follows:
We define the set $V_d=V_d(m)$ of vertices as 
\[V_d(m):=\{{\rm vt}(m-s;a^{(s)}_1,a^{(s)}_2,\cdots,a^{(s)}_d)|0\leq s\leq m,\
 a^{(s)}_i\in J \}. 
\]
And we define the set $E_d=E_d(m)$ of directed edges as 
\begin{multline*} E_d(m):=\{{\rm vt}(m-s;a^{(s)}_1,\cdots,a^{(s)}_d)\rightarrow
{\rm vt}(m-s-1;a^{(s+1)}_1,\cdots,a^{(s+1)}_d)\\
|\ 0\leq s\leq m-1,\ {\rm vt}(m-s;a^{(s)}_1,\cdots,a^{(s)}_d),\ {\rm vt}(m-s-1;a^{(s+1)}_1,\cdots,a^{(s+1)}_d)\in V_d(m)\}.
\end{multline*}
\end{defn}

Now, let us define the set of directed paths from ${\rm vt}(m;1,2,\cdots,d)$ to ${\rm vt}(0;m'+1,m'+2,\cdots,r,\ovl{d-r+m'},\ovl{d-r+m'-1},\cdots,\ovl{2},\ovl{1})$ (resp. ${\rm vt}(0;m'+1,m'+2,\cdots,m'+d)$) in the case $m'+d>r$ (resp. $m'+d\leq r$) in $(V_d,E_d)$.

\begin{defn}\label{pathdef}
Let $X_d(m,m')$ be the set of directed paths $p$
\begin{multline*}p={\rm vt}(m;a^{(0)}_1,\cdots,a^{(0)}_d)
\rightarrow{\rm vt}(m-1;a^{(1)}_1,\cdots,a^{(1)}_d)\rightarrow {\rm vt}(m-2;a^{(2)}_1,
\cdots,a^{(2)}_d)\\
\rightarrow
\cdots\rightarrow{\rm vt}(1;a^{(m-1)}_1,\cdots,a^{(m-1)}_d)
\rightarrow{\rm vt}(0;a^{(m)}_1,\cdots,a^{(m)}_d),
\end{multline*}
which satisfy the following conditions: For $0\leq s\leq m$,
\begin{enumerate}
\item $a^{(s)}_{\zeta}\in J$ $(1\leq \zeta\leq d)$, 
\item $a^{(s)}_{1}<a^{(s)}_{2}<\cdots<a^{(s)}_{d}$, 
\item If $a^{(s)}_{\zeta}\in\{j|1\leq j\leq r-1\}$, then $a^{(s+1)}_{\zeta}=a^{(s)}_{\zeta}$ or $a^{(s)}_{\zeta}+1$. If $a^{(s)}_{\zeta}=r$, then $a^{(s)}_{\zeta+1}\in\{r,\ovl{r},\ovl{r-1},\cdots,\ovl{1}\}$. If $a^{(s)}_{\zeta}\in\{\ovl{j}|1\leq j\leq r\}$, then $a^{(s+1)}_{\zeta}\in\{\ovl{|a^{(s)}_{\zeta}|},\ovl{|a^{(s)}_{\zeta}|-1},\cdots,\ovl{2},\ovl{1}\}$,
\item $(a^{(0)}_1,a^{(0)}_2,\cdots,a^{(0)}_d)=(1,2,\cdots,d)$,
\[ (a^{(m)}_1,\cdots,a^{(m)}_d)=
\begin{cases}
(m'+1,m'+2,\cdots,r,\ovl{d-r+m'},\cdots,\ovl{2},\ovl{1}) & {\rm if}\ m'+d>r,\\
(m'+1,m'+2,\cdots,m'+d) & {\rm if}\ m'+d\leq r,
\end{cases}
\]
\item If $a^{(s+1)}_{\zeta}\in\{\ovl{j}|\ 1\leq j\leq r\}$, then $|a^{(s+1)}_{\zeta}|>|a^{(s)}_{\zeta+1}|$.
\end{enumerate}
\end{defn}

\begin{defn}\label{connected-def}
We say that two vertices ${\rm vt}(m-s;a^{(s)}_1,\cdots,a^{(s)}_d)$ and ${\rm vt}(m-s-1;a^{(s+1)}_1,\cdots,a^{(s+1)}_d)$ are {\rm connected} if these vertices satisfy the conditions (i), (ii), (iii) and (v) in Definition \ref{pathdef}.
\end{defn}


Define a Laurent monomial associated with each edge of a path in $X_d(m,m')$.
\begin{defn}\label{labeldef}
Let $p\in X_d(m,m')$ be a path:
\begin{multline*}p={\rm vt}(m;a^{(0)}_1,\cdots,a^{(0)}_d)
\rightarrow{\rm vt}(m-1;a^{(1)}_1,\cdots,a^{(1)}_d)\rightarrow {\rm vt}(m-2;a^{(2)}_1,
\cdots,a^{(2)}_d)\\
\rightarrow
\cdots\rightarrow{\rm vt}(1;a^{(m-1)}_1,\cdots,a^{(m-1)}_d)
\rightarrow{\rm vt}(0;a^{(m)}_1,\cdots,a^{(m)}_d).
\end{multline*}
\begin{enumerate}
\item For each $0\leq s\leq m$, we define the {\it label of the edge} ${\rm vt}(m-s;a^{(s)}_1,a^{(s)}_2,\cdots,a^{(s)}_d)\rightarrow{\rm vt}(m-s-1;a^{(s+1)}_1,a^{(s+1)}_2,\cdots,a^{(s+1)}_d)$ as the Laurent monomial which determined as follows and denote it $Q^{(s)}(p)$ : We suppose that $0\leq \delta\leq d$, $1\leq a^{(s)}_1<\cdots<a^{(s)}_{\delta}\leq r$, and $a^{(s)}_{\delta+1},\cdots,a^{(s)}_d \in \{\ovl{i}|\ 1\leq i\leq r\}$. In the case $a^{(s)}_{\delta}<r$, we set
\[
Q^{(s)}(p):=\frac{Y_{m-s,a^{(s+1)}_1-1}}{Y_{m-s,a^{(s)}_1}}\cdots \frac{Y_{m-s,a^{(s+1)}_{\delta}-1}}{Y_{m-s,a^{(s)}_{\delta}}} \frac{Y_{m-s,|a^{(s)}_{\delta+1}|}}{Y_{m-s,|a^{(s+1)}_{\delta+1}|-1}} \cdots \frac{Y_{m-s,|a^{(s)}_{d}|}}{Y_{m-s,|a^{(s+1)}_d|-1}}.
\]
In the case $a^{(s)}_{\delta}=r$, we set
\[ Y(a^{(s+1)}_{\delta}):=
\begin{cases}
\frac{Y_{m-s,r-1}}{Y_{m-s,r}} & {\rm if}\ a^{(s+1)}_{\delta}=r,\\
\frac{1}{Y_{m-s,|a^{(s+1)}_{\delta}|-1}} & {\rm if}\ a^{(s+1)}_{\delta}\in \{\ovl{i} |i=1,\cdots,r \},
\end{cases}
\]
and set
\[
Q^{(s)}(p):=\frac{Y_{m-s,a^{(s+1)}_1-1}}{Y_{m-s,a^{(s)}_1}}\cdots \frac{Y_{m-s,a^{(s+1)}_{\delta-1}-1}}{Y_{m-s,a^{(s)}_{\delta-1}}}Y(a^{(s+1)}_{\delta}) \frac{Y_{m-s,|a^{(s)}_{\delta+1}|}}{Y_{m-s,|a^{(s+1)}_{\delta+1}|-1}} \cdots \frac{Y_{m-s,|a^{(s)}_{d}|}}{Y_{m-s,|a^{(s+1)}_d|-1}}.
\]
\item And we define the {\it label} $Q(p)$ {\it of the path} $p$ as the product of them:
\begin{equation}\label{alllabels} Q(p):=\prod_{s=0}^{m-1}Q^{(s)}(p). 
\end{equation}
\item For a subpath $p'$
\begin{multline*}
p'={\rm vt}(m-s';a^{(s')}_1,\cdots,a^{(s')}_d)\rightarrow{\rm vt}(m-s'-1;a^{(s'+1)}_1,
\cdots,a^{(s'+1)}_d)\rightarrow \\
\cdots\rightarrow{\rm vt}(m-s'-1';a^{(s''-1)}_1,\cdots,
a^{(s''-1)}_d)
\rightarrow{\rm vt}(m-s'';a^{(s'')}_1,\cdots,
a^{(s'')}_d)
\end{multline*}
of $p$ $(0\leq s'<s''\leq m)$, we define the {\it label} {\it of the subpath} $p'$ as 
\begin{equation}\label{sublabels} 
Q(p'):=\prod_{s=s'}^{s''}Q^{(s)}(p). 
\end{equation}
\end{enumerate}
\end{defn}

\begin{ex}\label{pathex}
Let $r=m=3$, $m'=2$, $d=2$. We can describe the paths of $X_2(3,2)$ as follows. For simplicity, we denote vertices ${\rm vt}(*;*,*)$ by $(*;*,*):$
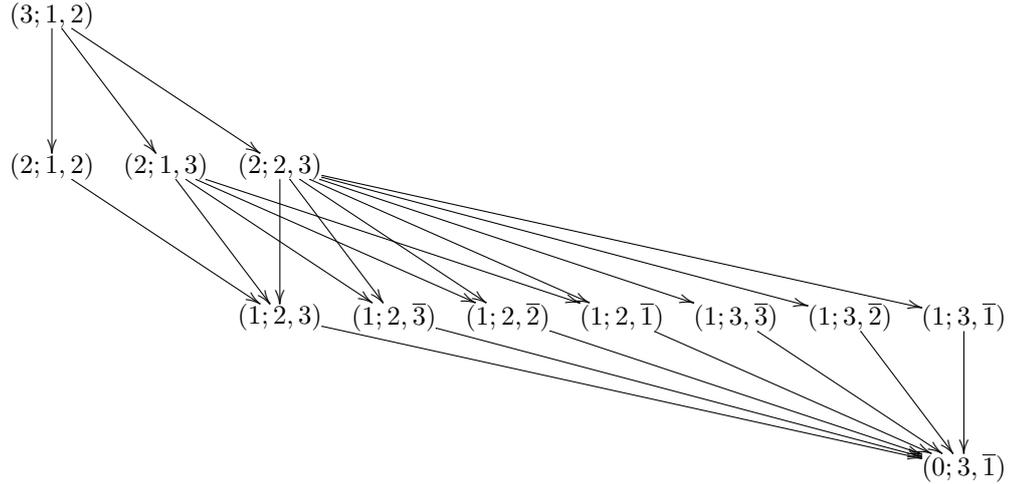
\[
\begin{xy}
(0,90) *{(3;1,2)}="3;1,2",
(0,70)*{(2;1,2)}="2;1,2",
(15,70)*{(2;1,3)}="2;1,3",
(30,70)*{(2;2,3)}="2;2,3",
(30,50)*{(1;2,3)}="1;2,3",
(45,50)*{(1;2,\ovl{3})}="1;2,4",
(60,50)*{(1;2,\ovl{2})}="1;2,5",
(75,50)*{(1;2,\ovl{1})}="1;2,6",
(90,50)*{(1;3,\ovl{3})}="1;3,4",
(105,50)*{(1;3,\ovl{2})}="1;3,5",
(120,50)*{(1;3,\ovl{1})}="1;3,6",
(120,30)*{(0;3,\ovl{1})}="0;3,4",
\ar@{->} "3;1,2";"2;1,2"
\ar@{->} "3;1,2";"2;1,3"
\ar@{->} "3;1,2";"2;2,3"
\ar@{->} "2;1,2";"1;2,3"
\ar@{->} "2;1,3";"1;2,3"
\ar@{->} "2;1,3";"1;2,4"
\ar@{->} "2;1,3";"1;2,5"
\ar@{->} "2;1,3";"1;2,6"
\ar@{->} "2;2,3";"1;2,3"
\ar@{->} "2;2,3";"1;2,4"
\ar@{->} "2;2,3";"1;2,5"
\ar@{->} "2;2,3";"1;2,6"
\ar@{->} "2;2,3";"1;3,4"
\ar@{->} "2;2,3";"1;3,5"
\ar@{->} "2;2,3";"1;3,6"
\ar@{->} "1;2,3";"0;3,4"
\ar@{->} "1;2,4";"0;3,4"
\ar@{->} "1;2,5";"0;3,4"
\ar@{->} "1;2,6";"0;3,4"
\ar@{->} "1;3,4";"0;3,4"
\ar@{->} "1;3,5";"0;3,4"
\ar@{->} "1;3,6";"0;3,4"
\end{xy}
\]
Thus, $X_2(3,2)$ has the following paths$:$

$p_1=(3;1,2)\rightarrow(2;1,2)\rightarrow(1;2,3)\rightarrow(0;3,\ovl{1})$,

$p_2=(3;1,2)\rightarrow(2;1,3)\rightarrow(1;2,3)\rightarrow(0;3,\ovl{1})$,

$p_3=(3;1,2)\rightarrow(2;1,3)\rightarrow(1;2,\ovl{3})\rightarrow(0;3,\ovl{1})$,

$p_4=(3;1,2)\rightarrow(2;1,3)\rightarrow(1;2,\ovl{2})\rightarrow(0;3,\ovl{1})$,

$p_5=(3;1,2)\rightarrow(2;1,3)\rightarrow(1;2,\ovl{1})\rightarrow(0;3,\ovl{1})$,
 
$p_6=(3;1,2)\rightarrow(2;2,3)\rightarrow(1;2,3)\rightarrow(0;3,\ovl{1})$,

$p_7=(3;1,2)\rightarrow(2;2,3)\rightarrow(1;2,\ovl{3})\rightarrow(0;3,\ovl{1})$,

$p_8=(3;1,2)\rightarrow(2;2,3)\rightarrow(1;2,\ovl{2})\rightarrow(0;3,\ovl{1})$,

$p_9=(3;1,2)\rightarrow(2;2,3)\rightarrow(1;2,\ovl{1})\rightarrow(0;3,\ovl{1})$,

$p_{10}=(3;1,2)\rightarrow(2;2,3)\rightarrow(1;3,\ovl{3})\rightarrow(0;3,\ovl{1})$,

$p_{11}=(3;1,2)\rightarrow(2;2,3)\rightarrow(1;3,\ovl{2})\rightarrow(0;3,\ovl{1})$,

$p_{12}=(3;1,2)\rightarrow(2;2,3)\rightarrow(1;3,\ovl{1})\rightarrow(0;3,\ovl{1})$.

\vspace{3mm}

Let us calculate the label of the path $p_1$. 
By Definition \ref{labeldef} $(iii)$, the label $Q^{(0)}(p_1)$ of the edge $(3;1,2)\rightarrow(2;1,2)$ is
\[ Q^{(0)}(p_1)=\frac{Y_{3,1-1}}{Y_{3,1}}\frac{Y_{3,2-1}}{Y_{3,2}}=\frac{1}{Y_{3,2}}, \] 
where we set $Y_{3,0}=1$ following Remark \ref{importantrem}. The labels of the edges $(2;1,2)\rightarrow(1;2,3)$ and $(1;2,3)\rightarrow(1;3,\ovl{1})$ are as follows:
\[ Q^{(1)}(p_1)=\frac{Y_{2,2-1}}{Y_{2,1}}\frac{Y_{2,3-1}}{Y_{2,2}}=1,\ \ 
Q^{(2)}(p_1)=\frac{Y_{1,3-1}}{Y_{1,2}}\frac{1}{Y_{1,1-1}}=1. \]
Therefore, we get $Q(p_1)=\frac{1}{Y_{3,2}}$.

Similarly, we have
\[ Q(p_1)=\frac{1}{Y_{3,2}},\ \ Q(p_2)=\frac{Y_{2,2}}{Y_{3,1}Y_{2,3}},\ 
\ Q(p_3)=\frac{Y_{1,3}}{Y_{3,1}Y_{2,2}},\ \ Q(p_4)=\frac{Y_{1,2}}{Y_{3,1}Y_{2,1}}, \]
\[ Q(p_5)=\frac{Y_{1,1}}{Y_{3,1}},\ \ Q(p_6)=\frac{Y_{2,1}}{Y_{2,3}},
\ \ Q(p_7)=\frac{Y_{2,1}Y_{1,3}}{Y_{2,2}^2},\ \ Q(p_8)=\frac{Y_{1,2}}{Y_{2,2}}, \]
\[ Q(p_9)=\frac{Y_{2,1}Y_{1,1}}{Y_{2,2}},\ \ Q(p_{10})=\frac{Y_{1,2}}{Y_{2,2}},
\ \ Q(p_{11})=\frac{Y_{1,2}^2}{Y_{2,1}Y_{1,3}},\ \ Q(p_{12})=\frac{Y_{1,1}Y_{1,2}}{Y_{1,3}}. \]

\end{ex}

\begin{defn}\label{iseq}
For each path $p\in X_d(m,m')$
\begin{multline*}p={\rm vt}(m;a^{(0)}_1,\cdots,a^{(0)}_d)\rightarrow{\rm vt}(m-1;a^{(1)}_1,\cdots,a^{(1)}_d)\rightarrow{\rm vt}(m-2;a^{(2)}_1,\cdots,a^{(2)}_d)\\
\rightarrow
\cdots\rightarrow{\rm vt}(1;a^{(m-1)}_1,\cdots,a^{(m-1)}_d)\rightarrow{\rm vt}(0;a^{(m)}_1,\cdots,a^{(m)}_d)
\end{multline*}
and $i\in \{1,\cdots,d\}$, we call the following sequence
\[ a^{(0)}_i\rightarrow a^{(1)}_i\rightarrow a^{(2)}_i\rightarrow\cdots\rightarrow a^{(m)}_i \]
an {\it i-sequence} of $p$.
\end{defn}

We can easily see the following by Definition \ref{pathdef} (iii) and (iv): For $1\leq i\leq d$,
\begin{equation}\label{iseqine}
i=a^{(0)}_i\leq a^{(1)}_i\leq \cdots\leq a^{(m)}_i,
\end{equation}
in the order (\ref{order}).

\subsection{One-to-one correspondence between paths in $X_d(m,m')$ and terms of $\Delta^L(k;\textbf{i})(\textbf{Y})$}

In this section, we describe the terms in $\Delta^L(k;\textbf{i})(\textbf{Y})$ as the paths in $X_d(m,m')$:

\begin{prop}\label{pathlem}We use the setting and the notations in Sect.\ref{gmc}:
\[ u=(s_1s_2\cdots s_r)^{m-1}s_1\cdots s_{i_n},\q v=e.\]
Then, we have the following:
\[ \Delta^L(k;\textbf{i})(\textbf{Y})=\sum_{p\in X_d(m,m')} Q(p). \]
\end{prop}

Let us give an overview of the proof of Proposition \ref{pathlem}. For $1\leq s\leq m$, we define
\begin{equation}\label{xmdef1}
x^{(s)}_{-[1,r]}:=x_{-1}(Y_{s,1})\cdots x_{-r}(Y_{s,r}).
\end{equation}
For $1\leq s\leq m$ and $i_1,\cdots,i_d\in J:=\{i,\ovl i|1\leq i\leq r\}$, we set 
\begin{equation}\label{xmdef2}
 (s;i_1,i_2,\cdots,i_d):=\lan x^{(1)}_{-[1,r]} x^{(2)}_{-[1,r]}\cdots x^{(s)}_{-[1,r]}
(v_{i_1}\wedge\cdots\wedge v_{i_d}), \q u_{\leq k}(v_1\wedge\cdots\wedge v_d)  \ran . 
\end{equation}

We shall prove $\Delta^L(k;\textbf{i})(\textbf{Y})=(m;1,2,\cdots,d)$ in Lemma \ref{xmlem1} (i). In Lemma \ref{xmlem1} (ii) and (iii), we shall also prove the recurrence formula for $\{(s;i_1,\cdots,i_d)\}$, which implies that $\Delta^L(k;\textbf{i})(\textbf{Y})=(m;1,2,\cdots,d)$ is expressed as a linear combination of
$\{(0;j_1,\cdots,j_d)| j_1,\cdots,j_d\in J,\ j_1<\cdots<j_d\}$. Note that if $(j_1,\cdots,j_d)=(m'+1,m'+2,\cdots,r,\ovl{d-r+m'},\ovl{d-r+m'-1},\cdots\ovl{1})$ (resp. $=(m'+1,m'+2,\cdots,m'+d)$), then $(0;j_1,\cdots,j_d)=1$ in the case $m'+d>r$ (resp. $m'+d\leq r$) by (\ref{ukeq1}), (\ref{ukeq2}) and (\ref{xmdef2}). If $(j_1,\cdots,j_d)$ is not as above, then we get $(0;j_1,\cdots,j_d)=0$. As a sequence of this calculation, we obtain Proposition \ref{pathlem} .

First, let us see the following lemma. We can verify it in the same way as (\ref{xmpro}).

\begin{lem}\label{xmpro2}
\[
x_{-i}(Y)v_j=
\begin{cases}
Y^{-1}v_{i}+v_{i+1} & {\rm if}\ j=i, \\
Y v_{i+1} & {\rm if}\ j=i+1, \\
v_j & {\rm otherwise},
\end{cases}\q
x_{-i}(Y)v_{\ovl{j}}=
\begin{cases}
Y^{-1}v_{\ovl{i+1}}+v_{\ovl{i}} & {\rm if}\ j=i+1, \\
Y v_{\ovl{i}} & {\rm if}\ j=i, \\
v_{\ovl{j}} & {\rm otherwise},
\end{cases}
\]
for all $1\leq i,j\leq r$ and $Y\in \mathbb{C}^{\times}$, where we set $v_{r+1}:=v_{\ovl{r}}$. 
\end{lem}

In the next lemma, we set $|l|=|\ovl{l}|=l$ for $1\leq l\leq r$.

\begin{lem}\label{xmlem1}
\begin{enumerate}
\item $\Delta^L(k;\textbf{i})(\textbf{Y})=(m;1,\cdots,d)$.

\item For $0\leq \delta\leq d$, $1\leq i_1<\cdots<i_{\delta}\leq r,\ i_{\delta+1},\cdots,i_d \in \{\ovl{i}|\ 1\leq i\leq r\}$ and $1\leq s\leq m$, we have the followings:

In the case $i_{\delta}<r$, 
\begin{eqnarray}\label{xmlem12}
& &(s;i_1,\cdots,i_{\delta},i_{\delta+1},\cdots,i_d) \nonumber \\
&=&\sum_{(j_1,\cdots,j_d)\in V}\frac{Y_{s,j_1-1}}{Y_{s,i_1}}\cdots \frac{Y_{s,j_{\delta}-1}}{Y_{s,i_{\delta}}} \frac{Y_{s,|i_{\delta+1}|}}{Y_{s,|j_{\delta+1}|-1}} \cdots \frac{Y_{s,|i_{d}|}}{Y_{s,|j_{d}|-1}} \qq \qq
\\
& &\qq \qq \qq \qq \qq \qq\cdot(s-1;j_1,\cdots,j_{\delta},j_{\delta+1},\cdots,j_d),  \nonumber
\end{eqnarray}
where $(j_1,\cdots,j_d)$ runs over $V:=\{(j_1,\cdots,j_d) |\ j_1<\cdots<j_{\delta},\ j_{\zeta}=i_{\zeta}\ {\rm or}\ i_{\zeta}+1\ (1\leq \zeta\leq \delta),\ j_{\zeta}\in\{\ovl{|i_{\zeta}|},\ \ovl{|i_{\zeta}|-1},\cdots,\ovl{1}\}\ (\delta+1\leq \zeta\leq d)\}$. 

In the case $i_{\delta}=r$, we set
\[ Y(j_{\delta}):=
\begin{cases}
\frac{Y_{s,r-1}}{Y_{s,r}} & {\rm if}\ j_{\delta}=r,\\
\frac{1}{Y_{s,|j_{\delta}|-1}} & {\rm if}\ j_{\delta}\in \{\ovl{i} |i=1,\cdots,r \}.
\end{cases}
\]
Then we have
\begin{eqnarray}\label{xmlem13}
& &(s;i_1,\cdots,i_{\delta-1},r,i_{\delta+1},\cdots,i_d)\nonumber \\
&=&\sum_{(j_1,\cdots,j_d)\in V}\frac{Y_{s,j_1-1}}{Y_{s,i_1}}\cdots \frac{Y_{s,j_{\delta-1}-1}}{Y_{s,i_{\delta-1}}}\cdot Y(j_{\delta})\cdot
 \frac{Y_{s,|i_{\delta+1}|}}{Y_{s,|j_{\delta+1}|-1}} \cdots \frac{Y_{s,|i_{d}|}}{Y_{s,|j_{d}|-1}}\qq \qq \\
& &\qq \qq \qq \qq \qq \qq\cdot(s-1;j_1,\cdots,j_{\delta-1},j_{\delta},j_{\delta+1},\cdots,j_d), \nonumber
\end{eqnarray}
where $(j_1,\cdots,j_d)$ runs over $V:=\{(j_1,\cdots,j_d) |\ j_1<\cdots<j_{\delta},\ j_{\zeta}=i_{\zeta}\ {\rm or}\ i_{\zeta}+1\ (1\leq \zeta\leq \delta-1),\ 
j_{\delta}\in \{r,\ \ovl{r},\ \ovl{r-1},\cdots,\ \ovl{1}\},\   
j_{\zeta}\in\{\ovl{|i_{\zeta}|},\ \ovl{|i_{\zeta}|-1},\cdots,\ovl{1}\}\ (\delta+1\leq \zeta\leq d)\}$. 

\item In addition to the assumptions in $(ii)$, we suppose that $i_1<\cdots<i_{\delta}<i_{\delta+1}<\cdots<i_d$ with the order $(\ref{order})$. If $i_{\delta}<r$, then we can reduce the range $V$ of the sum in $(\ref{xmlem12})$ to
\[ V':=\{(j_1,\cdots,j_d)\in V|\ |j_{l}|>|i_{l+1}|\ (\delta+1\leq l\leq d-1) \}. \] 
If $i_{\delta}=r$, then we can reduce the range $V$ of the sum in $(\ref{xmlem13})$ to
\[V':=
\begin{cases}
\{(j_1,\cdots,j_d)\in V|\ |j_{l}|>|i_{l+1}|\ (\delta+1\leq l\leq d-1) \}\ &\ {\rm if}\ j_{\delta}=r,\\ \{(j_1,\cdots,j_d)\in V|\ |j_{l}|>|i_{l+1}|\ (\delta\leq l\leq d-1) \}\ & \ {\rm if}\ j_{\delta}\in \{\ovl{i}|1\leq i\leq r\}. 
\end{cases}
\]
 
\end{enumerate}

\end{lem}
\nd
{\sl Proof.}

(i) By Lemma \ref{xmpro2}, if $i>j$ $(i,j\in\{1,\cdots,r\})$, then we have $x_{-i}(Y)v_j=v_j$. Thus, we get
\begin{eqnarray*}
& &(m;1,\cdots,d):=\lan x^{(1)}_{-[1,r]} \cdots x^{(m-1)}_{-[1,r]}x^{(m)}_{-[1,r]}(v_1\wedge\cdots \wedge v_d),\ u_{\leq k}(v_1\wedge\cdots\wedge v_d)\ran \\
&=&\lan x^{(1)}_{-[1,r]} \cdots x^{(m-1)}_{-[1,r]}x_{-1}(Y_{m,1})\cdots x_{-d}(Y_{m,d})(v_1\wedge\cdots \wedge v_d),\ u_{\leq k}(v_1\wedge\cdots\wedge v_d)\ran\\
 &=&\lan x^L_{\textbf{i}}(\textbf{Y})(v_1\wedge\cdots \wedge v_d),\q u_{\leq k}(v_1\wedge\cdots\wedge v_d)\ran=\Delta^L(k;\textbf{i})(\textbf{Y}). 
\end{eqnarray*}

(ii) By Lemma \ref{xmpro2}, for $1\leq s\leq m$ and $1\leq i\leq r$, we get
\begin{equation}\label{xmpro3}
x^{(s)}_{-[1,r]}v_i=
\begin{cases}
\frac{Y_{s,i-1}}{Y_{s,i}}v_i+v_{i+1} & {\rm if}\ 1\leq i\leq r-1, \\
\frac{Y_{s,r-1}}{Y_{s,r}}v_r+\sum^r_{j=1}\frac{1}{Y_{s,j-1}}v_{\ovl{j}}& {\rm if}\ i=r,
\end{cases}
\end{equation}
and
\begin{equation}\label{xmpro4}
x^{(s)}_{-[1,r]}v_{\ovl{i}}=\sum^i_{j=1}\frac{Y_{s,i}}{Y_{s,j-1}}v_{\ovl{j}},
\end{equation}
where we set $Y_{s,0}=1$. Since we supposed that $i_{\delta+1},\cdots,i_d \in \{\ovl{i}|\ 1\leq i\leq r\}$, if $i_{\delta}<r$, then
\begin{eqnarray}\label{xmpro6} 
& &x^{(1)}_{-[1,r]} \cdots x^{(s-1)}_{-[1,r]}x^{(s)}_{-[1,r]}
(v_{i_1}\wedge\cdots \wedge v_{i_{\delta}} \wedge v_{i_{\delta+1}}\wedge\cdots \wedge v_{i_d}) \nonumber \\
&=&x^{(1)}_{-[1,r]} \cdots x^{(s-1)}_{-[1,r]}(\left(\frac{Y_{s,i_1-1}}{Y_{s,i_1}}v_{i_1}+v_{i_1+1}\right)\wedge \cdots \wedge \left(\frac{Y_{s,i_{\delta}-1}}{Y_{s,i_{\delta}}}v_{i_{\delta}}+v_{i_{\delta}+1}\right) \nonumber\\
& &\qq \qq \qq \qq \qq \wedge \left(\sum^{|i_{\delta+1}|}_{l=1}\frac{Y_{s,|i_{\delta+1}|}}{Y_{s,l-1}}v_{\ovl{l}}\right) \wedge \cdots
\wedge \left(\sum^{|i_{d}|}_{l=1}\frac{Y_{s,|i_{d}|}}{Y_{s,l-1}}v_{\ovl{l}}\right)) \nonumber\\
&=&\sum_{j_1,\cdots,j_d} (\frac{Y_{s,j_1-1}}{Y_{s,i_1}}\cdots \frac{Y_{s,j_{\delta}-1}}{Y_{s,i_{\delta}}} \frac{Y_{s,|i_{\delta+1}|}}{Y_{s,|j_{\delta+1}|-1}} \cdots \frac{Y_{s,|i_{d}|}}{Y_{s,|j_{d}|-1}} \nonumber\\
& &\qq \qq x^{(1)}_{-[1,r]} \cdots x^{(s-1)}_{-[1,r]} (v_{j_1}\wedge\cdots\wedge v_{j_{\delta}}\wedge v_{j_{\delta+1}}\wedge\cdots\wedge v_{j_d})),
\end{eqnarray}
where $(j_1,\cdots,j_d)$ runs over $\{(j_1,\cdots,j_d) |\ j_1<\cdots<j_{\delta},\ j_{\zeta}=i_{\zeta}\ {\rm or}\ i_{\zeta}+1\ (1\leq \zeta\leq \delta),\ j_{\zeta}\in\{\ovl{|i_{\zeta}|},\ \ovl{|i_{\zeta}|-1},\cdots,\ovl{1}\}\ (\delta+1\leq \zeta\leq d)\}$. We remark that
\[ \frac{Y_{s,j_{\zeta}-1}}{Y_{s,i_{\zeta}}}=
\begin{cases}
\frac{Y_{s,i_{\zeta}-1}}{Y_{s,i_{\zeta}}} & {\rm if}\ j_{\zeta}=i_{\zeta},\\
1 & {\rm if}\ j_{\zeta}=i_{\zeta}+1,
\end{cases}
\]
for $1\leq \zeta\leq\delta$.

By pairing both sides in (\ref{xmpro6}) with $u_{\leq k}(v_1\wedge\cdots\wedge v_d)$, we obtain (\ref{xmlem12}). Similarly, we see (\ref{xmlem13}) in the case $i_{\delta}=r$. 

(iii) We suppose that $i_{\delta}<r$. Let $\hat{V}:=V\setminus V'$ be the complementary set. We define the map $\tau:\hat{V}\rightarrow \hat{V}$ as follows: Take $(j_1,\cdots,j_{\delta},j_{\delta+1},\cdots,j_d)\in \hat{V}$. Let $l$ $(\delta+1\leq l\leq d-1)$ be the index such that $|j_{\delta+1}|>|i_{\delta+2}|, \cdots, |j_{l-1}|>|i_{l}|$ and $|j_{l}|\leq |i_{l+1}|$. Since $|j_{l+1}|\leq |i_{l+1}|$ by the definition of $V$, we have $(j_1,\cdots,j_{l+1},j_{l},\cdots,j_d)\in \hat{V}$. So, we define $\tau(j_1,\cdots,j_{l},j_{l+1},\cdots,j_d):= (j_1,\cdots,j_{l+1},j_{l},\cdots,j_d)$. We can easily see that $\tau^2=id_{\hat{V}}$.

In $(\ref{xmlem12})$, $(s-1,j_1,\cdots,j_l,j_{l+1},\cdots,j_d)$ and $(s-1,j_1,\cdots,j_{l+1},j_l,\cdots,j_d)$ have the same coefficient
\[ \frac{Y_{s,j_1-1}}{Y_{s,i_1}}\cdots \frac{Y_{s,|i_l|}}{Y_{s,|j_{l}|-1}} \frac{Y_{s,|i_{l+1}|}}{Y_{s,|j_{l+1}|-1}} \cdots \frac{Y_{s,|i_{d}|}}{Y_{s,|j_{d}|-1}}=\frac{Y_{s,j_1-1}}{Y_{s,i_1}}\cdots  \frac{Y_{s,|i_{l}|}}{Y_{s,|j_{l+1}|-1}}\frac{Y_{s,|i_{l+1}|}}{Y_{s,|j_{l}|-1}} \cdots \frac{Y_{s,|i_{d}|}}{Y_{s,|j_{d}|-1}}. \]
Furthermore, by (\ref{xmdef2}), we obtain 
\[(s-1,j_1,\cdots,j_l,j_{l+1},\cdots,j_d)=-(s-1,j_1,\cdots,j_{l+1},j_l,\cdots,j_d). \]
Therefore, we get $\Sigma_{\hat{V}}=0$ in $(\ref{xmlem12})$, which implies our desired result. We can verify the case $i_{\delta}=r$ in the same way.
\qed

\vspace{3mm}

\nd
{\sl Proof of Proposition \ref{pathlem}.}

By the definition of $V$ and $V'$ in Lemma \ref{xmlem1}, we see that $(j_1,\cdots,j_d)\in V'$ if and only if the vertices vt$(s-1;j_1,\cdots,j_d)$ and vt$(s;i_1,\cdots,i_d)$ are connected (Definition \ref{connected-def}). Further, the coefficient of $(s-1;j_1,\cdots,j_d)$ in $(\ref{xmlem12})$, $(\ref{xmlem13})$ coincides with the label of the edge between vt$(s;i_1,\cdots,i_d)$ and vt$(s-1;j_1,\cdots,j_d)$ (Definition \ref{labeldef} (i)). Let us denote it by $^{(s)}Q^{i_1,\cdots,i_d}_{j_1,\cdots,j_d}$. Hence, in the case both $i_{\delta}=r$ and $i_{\delta}<r$, we get
\begin{equation}\label{plpr1}
 (s;i_1,\cdots,i_d)=\sum_{(j_1,\cdots,j_d)} \ ^{(s)}Q^{i_1,\cdots,i_d}_{j_1,\cdots,j_d}\ \cdot(s-1;j_1,\cdots,j_d),
\end{equation}
where $(j_1,\cdots,j_d)$ runs over the set $\{(j_1,\cdots,j_d)|\ {\rm vt}(s-1;j_1,\cdots,j_d)$ and ${\rm vt}(s;i_1$,
$\cdots,i_d)\ {\rm are\ connected}\}$. Note that the conditions $|j_l|>|i_{l+1}|$ in $V'$ and $|i_{l+1}|\geq |j_{l+1}|$ in $V$ implies $|j_l|>|j_{l+1}|$, and we get $j_1<j_2<\cdots<j_d$. Using Lemma \ref{xmlem1} (iii), we obtain the followings in the same way as (\ref{plpr1}):
\begin{equation}\label{plpr2}
 (s-1;j_1,\cdots,j_d)=\sum_{(k_1,\cdots,k_d)} \ ^{(s-1)}Q^{j_1,\cdots,j_d}_{k_1,\cdots,k_d}\ \cdot(s-2;k_1,\cdots,k_d),
\end{equation}
where $(k_1,\cdots,k_d)$ runs over the set $\{(k_1,\cdots,k_d)|\ {\rm vt}(s-2;k_1,\cdots,k_d)$ and ${\rm vt}(s-1;j_1,\cdots,j_d)\ {\rm are\ connected}\}$, and $^{(s-1)}Q^{j_1,\cdots,j_d}_{k_1,\cdots,k_d}$ is the label of the edge between vt$(s-1;j_1,\cdots,j_d)$ and vt$(s-2;k_1,\cdots,k_d)$. By (\ref{plpr1}), (\ref{plpr2}), $(s;i_1,\cdots,i_d)$ is a linear combination of $\{(s-2;k_1,\cdots,k_d)\}$, and the coefficient of $(s-2;k_1,\cdots,k_d)$ is as follows:
\[ \sum_{(j_1,\cdots,j_d)} \ ^{(s)}Q^{i_1,\cdots,i_d}_{j_1,\cdots,j_d}\cdot ^{(s-1)}Q^{j_1,\cdots,j_d}_{k_1,\cdots,k_d}\ \cdot(s-2;k_1,\cdots,k_d),\]
where $(j_1,\cdots,j_d)$ runs over the set $\{(j_1,\cdots,j_d)|\ {\rm vt}(s-1;j_1,\cdots,j_d)$ is connected to the vertices vt$(s;i_1,\cdots,i_d)$ and vt$(s-2;k_1,\cdots,k_d)\}$. The coefficient $^{(s)}Q^{i_1,\cdots,i_d}_{j_1,\cdots,j_d}\cdot ^{(s-1)}Q^{j_1,\cdots,j_d}_{k_1,\cdots,k_d}$ coincides with the label of subpath (Definition \ref{labeldef} (iii))
\[{\rm vt}(s;i_1,\cdots,i_d)\rightarrow {\rm vt}(s-1;j_1,\cdots,j_d)\rightarrow{\rm vt}(s-2;k_1,\cdots,k_d).\]

Repeating this argument, we see that $(s;i_1,\cdots,i_d)$ is a linear combination of $\{(0;l_1,\cdots,l_d)\}$ $(1\leq l_1<\cdots<l_d\leq \ovl{1})$. The coefficient of $(0;l_1,\cdots,l_d)$ is equal to the sum of labels of all subpaths from vt$(s;i_1,\cdots,i_d)$ to vt$(0;l_1,\cdots,l_d)$. In the case $m'+d>r$ (resp. $m'+d\leq r$), for $1\leq l_1<\cdots<l_d\leq \ovl{1}$, if $(l_1,\cdots,l_d)=(m'+1,m'+2,\cdots,r,\ovl{d-r+m'},\cdots,\ovl{2},\ovl{1})$ (resp. $=(m'+1,m'+2,\cdots,m'+d)$), then we obtain $(0;l_1,\cdots,l_d)=1$ by (\ref{ukeq1}), (\ref{ukeq2}) and (\ref{xmdef2}). If $(l_1,\cdots,l_d)$ is not as above, we obtain $(0;l_1,\cdots,l_d)=0$. Therefore, we see that $(s;i_1,\cdots,i_d)$ is equal to the sum of labels of subpaths from vt$(s;i_1,\cdots,i_d)$ to vt$(0;m'+1,m'+2,\cdots,r,\ovl{d-r+m'},\cdots,\ovl{2},\ovl{1})$ (resp. vt$(0;m'+1,m'+2,\cdots,m'+d))$.

In particular, $\Delta^L(k;\textbf{i})(\textbf{Y})=(m;1,2,\cdots,d)$ is equal to the sum of labels of paths in $X_d(m,m')$, which means $\Delta^L(k;\textbf{i})(\textbf{Y})=\sum_{p\in
X_d(m,m')} Q(p)$. 
\qed

\begin{ex}\label{pathex2}
Let us assume the same setting as Example \ref{pathex3}, {\it i.e.},
$r=3$, $u=s_1s_2s_3s_1s_2s_3s_1s_2$, $v=e$, $k=5$, 
$\textbf{i}=(-1,-2,-3,-1,-2,-3,-1,-2)$, 
$m=3$, $m'=2$ and $d=2$.
Therefore,  by Example \ref{pathex}, we obtain
\begin{eqnarray*}
\Delta^L(5;\textbf{i})(\textbf{Y})&=&
\frac{1}{Y_{3,2}}+\frac{Y_{2,2}}{Y_{3,1}Y_{2,3}}+\frac{Y_{1,3}}{Y_{3,1}Y_{2,2}}+\frac{Y_{1,2}}{Y_{3,1}Y_{2,1}}+\frac{Y_{1,1}}{Y_{3,1}}+\frac{Y_{2,1}}{Y_{2,3}} \\
& &+\frac{Y_{2,1}Y_{1,3}}{Y_{2,2}^2}+2\frac{Y_{1,2}}{Y_{2,2}}+\frac{Y_{2,1}Y_{1,1}}{Y_{2,2}}+\frac{Y_{1,2}^2}{Y_{2,1}Y_{1,3}}+\frac{Y_{1,1}Y_{1,2}}{Y_{1,3}}. 
\end{eqnarray*}
We find that this just coincides with the explicit form of 
$\Delta^L(5;\textbf{i})(\textbf{Y})$ in Example \ref{pathex3}.
\end{ex}

\begin{rem}\label{coinArem}
We suppose that $m'+d\leq r$.
\begin{enumerate}
\item[$(1)$] Definition \ref{pathdef} shows that the set $X_d(m,m')$ is constituted by paths $p$
\begin{multline*}p={\rm vt}(m;a^{(0)}_1,\cdots,a^{(0)}_d)\rightarrow{\rm vt}(m-1;a^{(1)}_1,\cdots,a^{(1)}_d)\rightarrow\\
\cdots\rightarrow{\rm vt}(1;a^{(m-1)}_1,\cdots,a^{(m-1)}_d)\rightarrow{\rm vt}(0;a^{(m)}_1,\cdots,a^{(m)}_d)
\end{multline*}
which satisfy the following conditions: For $0\leq s\leq m$,
\begin{enumerate}
\item[$(i)$] $a^{(s)}_{\zeta}\in\{1,\cdots,r\}$ $(1\leq \zeta\leq d)$, 
\item[$(ii)$] $a^{(s)}_{1}<a^{(s)}_{2}<\cdots<a^{(s)}_{d}$, 
\item[$(iii)$] $a^{(s+1)}_{\zeta}=a^{(s)}_{\zeta}$ or $a^{(s)}_{\zeta}+1$. 
\item[$(iv)$] $(a^{(0)}_1,a^{(0)}_2,\cdots,a^{(0)}_d)=(1,2,\cdots,d)$, 

$(a^{(m)}_1,\cdots,a^{(m)}_d)=(m'+1,m'+2,\cdots,m'+d)$.
\end{enumerate}
\item[$(2)$] By Definition \ref{labeldef}, the label $Q^{(s)}(p)$ of the edge ${\rm vt}(m-s;a^{(s)}_1,a^{(s)}_2,\cdots,a^{(s)}_d)\rightarrow{\rm vt}(m-s-1;a^{(s+1)}_1,a^{(s+1)}_2,\cdots,a^{(s+1)}_d)$ is as follows:
\[
Q^{(s)}(p):=\frac{Y_{m-s,a^{(s+1)}_1-1}}{Y_{m-s,a^{(s)}_1}}\cdots \frac{Y_{m-s,a^{(s+1)}_d-1}}{Y_{m-s,a^{(s)}_d}}.
\]
\item[$(3)$] For $G_A=SL_{r+1}(\mathbb{C})$, let $B_A$ and $(B_-)_A$ be two opposite Borel subgroups in $G_A$, $N_A\subset B_A$ and $(N_-)_A\subset (B_-)_A$ their unipotent radicals, and $W_A$ be the Weyl group of $G_A$. We define a reduced double Bruhat cell as $L^{u,v}_A:= N_A\cdot u\cdot N_A \cap (B_-)_A\cdot v\cdot (B_-)_A$.
We set $u,v\in W_A$ and their reduced word $\textbf{i}_A$ as
\[ u=\underbrace{s_1\cdots s_r}_{1\ {\rm st\ cycle}}\underbrace{s_1\cdots s_{r-1}}_{2\ {\rm nd\ cycle}}\cdots \underbrace{s_1\cdots s_{i_n}}_{m\ {\rm th\ cycle}}, \q v=e, \]
\[ \textbf{i}_A=(\underbrace{1,\cdots r}_{1\ {\rm st\ cycle}},\underbrace{1\cdots (r-1)}_{2\ {\rm nd\ cycle}}\cdots\underbrace{1,\cdots ,i_n}_{m\ {\rm th\ cycle}}), \]
where $n=l(u)$ and $1\leq i_n\leq r-m+1$. Let $i_k$ be the $k$-th index of  $\textbf{i}_A$ from the left, and belong to $m'$-th cycle. Using Theorem \ref{fp2}, we can define $\Delta^{L_A}(k;\textbf{i}_A)(\textbf{Y}_A):=(\Delta(k;\textbf{i}_A)\circ
 x^{L_A}_{\textbf{i}_A})(\textbf{Y}_A)$ in the same way as Definition \ref{gendef}, where 
\[
 \textbf{Y}_A:=(Y_{1,1},Y_{1,2},\cdots,Y_{1,r},Y_{2,1},Y_{2,2},\cdots,Y_{2,r-1},\cdots,Y_{m,1},\cdots,Y_{m,i_n})\in  (\mathbb{C}^{\times})^{n},
\]
and the map $x^{L_A}_{\textbf{i}_A}:(\mathbb{C}^{\times})^{n} \overset{\sim}{\rightarrow} L^{u,v}_A$ is defined as in Theorem \ref{fp2}.

Then, we already had seen in \cite{KaN} that 
$\Delta^{L_A}(k;\textbf{i}_A)(\textbf{Y}_A)=\sum_{p\in X_d(m,m')}Q(p)$, 
where $X_d(m,m')$ and the label $Q=\prod^{m-1}_{s=0}Q^{(s)}(p)$ 
is the one we have seen in $(1)$ and $(2)$. 
Therefore, it follows from Proposition \ref{pathlem} that 
if $m'+d\leq r$, then $\Delta^L(k;\textbf{i})(\textbf{Y})$ 
coincides with $\Delta^{L_A}(k;\textbf{i}_A)(\textbf{Y}_A)$.
\end{enumerate}
\end{rem}

\subsection{The properties of paths in $X_d(m,m')$}

In this subsection, we shall see some lemmas on $X_d(m,m')$. By Remark \ref{coinArem}, we suppose that $m'+d>r$. We fix a path $p\in X_d(m,m')$
\begin{multline}\label{fixpath}
p={\rm vt}(m;a^{(0)}_1,\cdots,a^{(0)}_d)\rightarrow
\cdots\rightarrow{\rm vt}(2;a^{(m-2)}_1,\cdots,a^{(m-2)}_d)  \\
\rightarrow{\rm vt}(1;a^{(m-1)}_1,\cdots,a^{(m-1)}_d)\rightarrow{\rm vt}(0;a^{(m)}_1,\cdots,a^{(m)}_d).
\end{multline}

\begin{lem}\label{fstlem}
For $p\in X_d(m,m')$ in $(\ref{fixpath})$,
$i$ $(1\leq i\leq d-1)$ and $s$ $(1\leq s\leq m)$, if $a^{(s)}_i\in\{\ovl{j}|1\leq j\leq r\}$, then we have $a^{(s-1)}_{i+1}\in\{\ovl{j}|1\leq j\leq r\}$ and
\[ a^{(s)}_i<a^{(s-1)}_{i+1}. \]
\end{lem}
\nd
{\sl Proof.}

Using Definition \ref{pathdef} (iii) and the assumption $a^{(s)}_i\in\{\ovl{j}|1\leq j\leq r\}$, we obtain $a^{(s-1)}_{i}\in \{r,\ovl{r},\ovl{r-1},\cdots,\ovl{1}\}$. Therefore, we also get $a^{(s-1)}_{i+1}\in \{\ovl{r},\ovl{r-1},\cdots,\ovl{1}\}$ by Definition \ref{pathdef} (ii). Further, it follows from Definition \ref{pathdef} (v) that $a^{(s)}_i<a^{(s-1)}_{i+1}$. \qed

\begin{lem}\label{cancellem}
For $p\in X_d(m,m')$ in $(\ref{fixpath})$ and $i$ $(r-m'+1\leq i\leq d)$, we obtain 
\begin{equation}\label{fixedeq}
a^{(m)}_i=a^{(m-1)}_i=\cdots=a^{(m-i+r-m'+1)}_i=\ovl{d-i+1}.
\end{equation}
\end{lem}
\nd
{\sl Proof.}

By Definition \ref{pathdef} (iv), we get $a^{(m)}_{r-m'+1}=\ovl{d-r+m'}$, and by Lemma \ref{fstlem}, we also get $\ovl{d-r+m'}=a^{(m)}_{r-m'+1}<a^{(m-1)}_{r-m'+2}\leq \ovl{1}$. Using Lemma \ref{fstlem} repeatedly, we obtain $\ovl{d-r+m'}=a^{(m)}_{r-m'+1}<a^{(m-1)}_{r-m'+2}<a^{(m-2)}_{r-m'+3}<\cdots<a^{(m-d+r-m'+1)}_{d}\leq \ovl{1}$, which means
\[a^{(m-i+r-m'+1)}_{i}=\ovl{d-i+1}\ \ (r-m'+1\leq i\leq d). \]
It follows from (\ref{iseqine}) and Definition \ref{pathdef} (iv) that $\ovl{d-i+1}=a^{(m-i+r-m'+1)}_{i}\leq a^{(m-i+r-m'+2)}_{i}\leq\cdots\leq a^{(m-1)}_{i}\leq a^{(m)}_{i}=\ovl{d-i+1}$, which yields (\ref{fixedeq}). \qed

\vspace{3mm}

By this lemma, we get $a^{(s)}_i=\ovl{d-i+1}$ for $r-m'+1\leq i$ and $m-i+r-m'+1\leq s\leq m$. In the next lemma, we see the properties for $a^{(s)}_i$ $(0\leq s\leq m-i+r-m')$.

\begin{lem}\label{tab-length}
For $i$ $(1\leq i\leq d)$ and $p\in X_d(m,m')$, let
\[ a^{(0)}_{i}\rightarrow a^{(1)}_{i}\rightarrow a^{(2)}_{i}\rightarrow \cdots\rightarrow a^{(m)}_{i} \]
be the $i$-sequence of the path $p$ $(Definition\ \ref{iseq})$.
\begin{enumerate}
\item In the case $i\leq r-m'$,
\[   \#\{0\leq s\leq m-1|\ 1\leq a^{(s)}_i\leq r,\ {\rm and}\ a^{(s)}_i=a^{(s+1)}_i\}=m-m'. \]
\item In the case $i> r-m'$,
\begin{multline*}
\#\{0\leq s\leq m-i+r-m'|\ 1\leq a^{(s)}_i\leq r\, {\rm and}\ a^{(s)}_i=a^{(s+1)}_i\}+\\
\#\{0\leq s\leq m-i+r-m'|\ \ovl{r}\leq a^{(s)}_i\leq \ovl{1}\}=m-m'. 
\end{multline*}
\end{enumerate}
\end{lem}
\nd
{\sl Proof.}

(i) In the case $i\leq r-m'$, Definition \ref{pathdef} (iv) and (\ref{iseqine}) show that 
\begin{equation}\label{tab-length1}
 i=a^{(0)}_i\leq a^{(1)}_i\leq\cdots\leq a^{(m)}_i=m'+i,\ \ a^{(s+1)}_i=a^{(s)}_i\ {\rm or}\ a^{(s)}_i+1. 
\end{equation}
In particular, we get $1\leq a^{(s)}_i\leq r$ for $1\leq s\leq m$.
By (\ref{tab-length1}), we obtain
\[ \#\{0\leq s\leq m-1|\ a^{(s+1)}_i=a^{(s)}_i+1\}=m',
\]
which implies $\#\{0\leq s\leq m-1|\ a^{(s)}_i=a^{(s+1)}_i\}=m-m'$. 

(ii) In the case $i>r-m'$, by (\ref{iseqine}), we have
\[
 i=a^{(0)}_i\leq a^{(1)}_i\leq\cdots\leq a^{(m-i+r-m')}_i\leq \ovl{1}. 
\]
We suppose that 
\begin{equation}\label{tab-length2}
i= a^{(0)}_i\leq a^{(1)}_i\leq \cdots\leq a^{(l)}_i\leq r,\ {\rm and}\ \ \ovl{r}\leq a^{(l+1)}_i\leq \cdots\leq a^{(m-i+r-m')}_i\leq\ovl{1},
\end{equation}
for some $1\leq l\leq m-i+r-m'$. Definition \ref{pathdef} (iii) implies that $a^{(s+1)}_i=a^{(s)}_i$ or $a^{(s)}_i+1$ $(1\leq s\leq l-1)$ and $a^{(l)}_i=r$. Therefore, 
\[ i=a^{(0)}_i\leq a^{(1)}_i\leq\cdots\leq a^{(l)}_i=r,\ a^{(s+1)}_i=a^{(s)}_i\ {\rm or}\ a^{(s)}_i+1. \]
So we have $\#\{1\leq s\leq l-1 |\ a^{(s+1)}_i=a^{(s)}_i \}=l-(r-i)$ in the same way as (i).

On the other hand, by the assumption $\ovl{r}\leq a^{(l+1)}_i\leq\cdots\leq a^{(m-i+r-m')}_i\leq \ovl{1}$ in (\ref{tab-length2}), we clearly see that $\#\{l+1\leq s\leq m-i+r-m' |\ \ovl{r}\leq a^{(s)}_i\leq \ovl{1} \}=m-i+r-m'-l$. Hence, $\#\{1\leq s\leq l-1 |\ a^{(s+1)}_i=a^{(s)}_i \}+\#\{l+1\leq s\leq m-i+r-m' |\ \ovl{r}\leq a^{(s)}_i\leq \ovl{1} \}=l-(r-i)+m-i+r-m'-l=m-m'$. \qed

\vspace{3mm}

By this lemma, we define $l^{(s)}_i\in\{0,1,\cdots,m\}$ $(1\leq i\leq d,\ 1\leq s\leq m-m')$ for the path $p\in X_d(m,m')$ in $(\ref{fixpath})$ as follows:  For $i\leq r-m'$, we set $\{l^{(s)}_i\}_{1\leq s\leq m-m'}$ $(l^{(1)}_i<\cdots<l^{(m-m')}_i)$ as
\begin{equation}\label{ldef}
\{l^{(1)}_i,\ l^{(2)}_i,\cdots,l^{(m-m')}_i\}:=\{s|a^{(s)}_i=a^{(s+1)}_i,\ \ 0\leq s\leq m-1 \}.
\end{equation}
For $i>r-m'$, we set $\{l^{(s)}_i\}_{1\leq s\leq m-m'}$ $(l^{(1)}_i<\cdots<l^{(m-m')}_i)$ as
\begin{eqnarray}\label{ldef2}
& &\{l^{(1)}_i,\ l^{(2)}_i,\cdots,l^{(m-m')}_i\}\nonumber \\
&:=&\{s|\ 1\leq a^{(s)}_i\leq r,\ a^{(s)}_i=a^{(s+1)}_i,\ 0\leq s\leq m-i+r-m' \}\nonumber \\
&\cup&\ \{s|\ \ovl{r}\leq a^{(s)}_i\leq \ovl{1}, \ 0\leq s\leq m-i+r-m' \}. \qq \q
\end{eqnarray}
We also set $k^{(s)}_i\in \{j,\ovl{j}|\ 1\leq j\leq r\}$ $(1\leq i\leq d,\ 1\leq s\leq m-m')$ as
\begin{equation}\label{kdef}
k^{(s)}_i:=a^{(l^{(s)}_i)}_i.
\end{equation}
Using (\ref{iseqine}) and $l^{(1)}_i<\cdots<l^{(m-m')}_i$, we obtain
\begin{equation}\label{kpro1}
k^{(1)}_i\leq \cdots\leq k^{(m-m')}_i.
\end{equation}
For $1\leq i\leq d$, let us define $\delta_i$ $(0\leq\delta_i\leq m-m')$ as
\begin{equation}\label{deldef}
1\leq k^{(1)}_i\leq\cdots\leq k^{(\delta_i)}_i\leq r<\ovl{r}\leq k^{(\delta_i+1)}_i\leq\cdots\leq k^{(m-m')}_i\leq\ovl{1},
\end{equation}
which is uniquely determined from $\{k^{(s)}_i\}_{s=1,\cdots,m-m'}$.

\begin{lem}\label{kpro2}
\begin{enumerate}
\item For $1\leq i\leq d$,
\[ l^{(s)}_{i}=
\begin{cases}
k^{(s)}_{i}+s-i-1 & {\rm if}\  k^{(s)}_{i}\in\{j|1\leq j\leq r\}, \\
s-i+r & {\rm if}\ k^{(s)}_{i}\in\{\ovl{j}|1\leq j\leq r\}.
\end{cases}
\]
\item For $1\leq s\leq m-m'$ and $1\leq i\leq d-1$, if $k^{(s)}_{i}\in\{j|1\leq j\leq r\}$, then
\[ k^{(s)}_{i}<k^{(s)}_{i+1},\q l^{(s)}_i\leq l^{(s)}_{i+1}. \]
For $1\leq i\leq d-1$, if $k^{(s)}_{i}\in\{\ovl{j}|1\leq j\leq r\}$, then
\[ k^{(s)}_{i}<k^{(s)}_{i+1},\q l^{(s)}_{i}=l^{(s)}_{i+1}+1. \]
\end{enumerate}
\end{lem}
\nd
{\sl Proof.}

(i) We suppose that $k^{(s)}_{i}\in \{j|1\leq j\leq r\}$. The definition of $l^{(s)}_i$ in (\ref{ldef}) means that the path $p$ has the following $i$-sequence (Definition \ref{iseq}):
\begin{flushleft}
$a^{(0)}_i=i,\ a^{(1)}_i=i+1,\ a^{(2)}_i=i+2,\cdots,a^{(l^{(1)}_i)}_i=i+l^{(1)}_i$,
 
$ a^{(l^{(1)}_i+1)}_i=i+l^{(1)}_i,\ a^{(l^{(1)}_i+2)}_i=i+l^{(1)}_i+1,\cdots, a^{(l^{(2)}_i)}_i=i+l^{(2)}_i-1, $

$ a^{(l^{(2)}_i+1)}_i=i+l^{(2)}_i-1,\ a^{(l^{(2)}_i+2)}_i=i+l^{(2)}_i,\cdots, a^{(l^{(3)}_i)}_i=i+l^{(3)}_i-2, $
\end{flushleft}
\begin{equation}\label{jlist} \vdots \end{equation}
\[ a^{(l^{(s-1)}_i+1)}_i=i+l^{(s-1)}_i-s+2,\ a^{(l^{(s-1)}_i+2)}_i=i+l^{(s-1)}_i-s+3,\cdots, a^{(l^{(s)}_i)}_i=i+l^{(s)}_i-s+1,\]
$a^{(l^{(s)}_i+1)}_i=i+l^{(s)}_i-s+1,\ a^{(l^{(s)}_i+2)}_i=i+l^{(s)}_i-s+2,\cdots$.

\vspace{2mm}

Hence we have 
\begin{equation}\label{kpro2pr1}
 k^{(s)}_{i}=a^{(l^{(s)}_{i})}_i=i+l^{(s)}_{i}-s+1,
\end{equation}
which implies $l^{(s)}_{i}=k^{(s)}_{i}+s-i-1$.

Next, we suppose that $a^{(l^{(s)}_{i})}_i=k^{(s)}_{i}\in \{\ovl{j}|1\leq j\leq r\}$. Using (\ref{iseqine}), we get $a^{(l^{(s)}_{i})}_i\leq a^{(l^{(s)}_{i}+1)}_i\leq\cdots\leq a^{(m-i+r-m')}_i$ and $a^{(\zeta)}_i\in\{\ovl{j}|1\leq j\leq r\}$ $(l^{(s)}_{i}\leq \zeta \leq m-i+r-m')$. Thus, by the definition (\ref{ldef2}) of $l^{(s)}_i$, we obtain $l^{(m-m')}_i=m-i+r-m'$, $l^{(m-m'-1)}_i=m-i+r-m'-1,l^{(m-m'-2)}_i=m-i+r-m'-2,\cdots, l^{(\xi)}_i=\xi-i+r$ $(s\leq\xi\leq m-m')$. In particular, we get
\begin{equation}\label{kpro2pr10}
 l^{(s)}_i=s-i+r.
\end{equation}  

(ii) We suppose that $k^{(s)}_{i}\in \{j|1\leq j\leq r\}$. If $k^{(s)}_{i+1}\in \{\ovl{j}|1\leq j\leq r\}$, then we obtain $k^{(s)}_{i}<k^{(s)}_{i+1}$ in the order (\ref{order}), and it follows $l^{(s)}_{i}\leq l^{(s)}_{i+1}$ from (i). So we may assume that $k^{(s)}_{i+1}\in \{j|1\leq j\leq r\}$. 

By Definition \ref{pathdef} (ii) and the definition (\ref{ldef2}) of $l^{(s)}_{i+1}$, we have $a^{(l^{(s)}_{i+1}+1)}_i<a^{(l^{(s)}_{i+1}+1)}_{i+1}=a^{(l^{(s)}_{i+1})}_{i+1}=k^{(s)}_{i+1}\leq r$. Therefore, the inequality (\ref{iseqine}) implies
\begin{equation}\label{kpro2pr3ano}
i=a^{(0)}_{i}\leq a^{(1)}_{i}\leq\cdots\leq a^{(l^{(s)}_{i+1})}_{i}\leq a^{(l^{(s)}_{i+1}+1)}_{i}<r, 
\end{equation}
\[ a^{(\zeta)}_{i}=a^{(\zeta-1)}_{i}\ {\rm or}\ a^{(\zeta-1)}_{i}+1\ (1\leq\zeta\leq l^{(s)}_{i+1}+1). \]
We obtain
\begin{equation}\label{kpro2pr4ano}
l^{(s)}_{i+1}+1-s\geq\#\{\zeta| a^{(\zeta)}_{i}=a^{(\zeta-1)}_{i}+1,\ 1\leq\zeta\leq l^{(s)}_{i+1}+1\},
\end{equation}
otherwise, it follows from (\ref{kpro2pr3ano}) and (i) that $a^{(l^{(s)}_{i+1}+1)}_i> i+l^{(s)}_{i+1}+1-s=k^{(s)}_{i+1}-1=a^{(l^{(s)}_{i+1})}_{i+1}-1$, and hence $a^{(l^{(s)}_{i+1}+1)}_i\geq a^{(l^{(s)}_{i+1})}_{i+1}$, which contradicts Definition \ref{pathdef} (ii).

The inequality (\ref{kpro2pr4ano}) means that
\begin{equation}\label{kpro2pr5ano}
s\leq\#\{\zeta| a^{(\zeta)}_{i}=a^{(\zeta-1)}_{i},\ 1\leq\zeta\leq l^{(s)}_{i+1}+1\}.
\end{equation}
On the other hand, the definition of $l^{(s)}_i$ implies $a^{(l^{(s)}_i+1)}_i=a^{(l^{(s)}_i)}_i=k^{(s)}_{i}\in\{j|1\leq j\leq r\}$. The inequality (\ref{iseqine}) shows 
\[
i=a^{(0)}_{i}\leq a^{(1)}_{i}\leq \cdots\leq a^{(l^{(s)}_{i})}_{i}=a^{(l^{(s)}_{i}+1)}_{i}=k^{(s)}_{i},
\] 
\[ a^{(\zeta)}_{i}=a^{(\zeta-1)}_{i}\ {\rm or}\ a^{(\zeta-1)}_{i}+1\ (1\leq\zeta\leq l^{(s)}_{i}+1),\]
and
\begin{equation}\label{kpro2pr6ano}
s=\#\{\zeta| a^{(\zeta)}_{i}=a^{(\zeta-1)}_{i},\ 1\leq\zeta\leq l^{(s)}_{i}+1\}.
\end{equation}
Since $a^{(l^{(s)}_i)}_i=a^{(l^{(s)}_i+1)}_i$, the equation (\ref{kpro2pr6ano}) means 
\begin{equation}\label{kpro2pr6ano2}
s-1=\#\{\zeta| a^{(\zeta)}_{i}=a^{(\zeta-1)}_{i},\ 1\leq\zeta\leq l^{(s)}_{i}\}.
\end{equation}
Thus, by (\ref{kpro2pr5ano}) and (\ref{kpro2pr6ano2}), we have $l^{(s)}_{i}< l^{(s)}_{i+1}+1$, and hence $l^{(s)}_{i}\leq l^{(s)}_{i+1}$, which yields $k^{(s)}_i<k^{(s)}_{i+1}$ since $k^{(s)}_i=i+l^{(s)}_i-s+1<i+l^{(s)}_{i+1}-s+2=(i+1)+l^{(s)}_{i+1}-s+1=k^{(s)}_{i+1}$.

Next, we suppose that $a^{(l^{(s)}_i)}_{i}=k^{(s)}_{i}\in \{\ovl{j}|1\leq j\leq r\}$. As we have seen in Lemma \ref{fstlem}, we obtain $a^{(l^{(s)}_i-1)}_{i+1}\in \{\ovl{j}|1\leq j\leq r\}$. Since $a^{(l^{(s)}_i-1)}_{i+1}\leq a^{(l^{(s)}_i)}_{i+1}\leq \cdots\leq a^{(m-(i+1)+r-m')}_{i+1}$, we get $a^{(\zeta)}_{i+1}\in\{\ovl{j}|1\leq j\leq r\}$ $(l^{(s)}_i-1\leq\zeta\leq m-(i+1)+r-m')$ and $l^{(m-m')}_{i+1}=m-(i+1)+r-m'$, $l^{(m-m'-1)}_{i+1}=m-(i+1)+r-m'-1,\cdots,l^{(\xi)}_{i+1}=\xi-(i+1)+r,\cdots$ $(s\leq\xi\leq m-m')$ by the definition (\ref{ldef2}) of $l^{(\xi)}_{i+1}$. In particular, we get
\[
l^{(s)}_{i+1}=s-(i+1)+r.
\]
Therefore, it follows from (\ref{kpro2pr10}) that $l^{(s)}_{i}=l^{(s)}_{i+1}+1$. Further, $k^{(s)}_i=a^{(l^{(s)}_i)}_i<a^{(l^{(s)}_i-1)}_{i+1}=a^{(l^{(s)}_{i+1})}_{i+1}=k^{(s)}_{i+1}$ by Lemma \ref{fstlem}. \qed

\subsection{The proof of Theorem \ref{thm1}}

In this subsection, we shall prove Theorem \ref{thm1}. First, we see the following lemma. Let us recall the definition (\ref{ccbar}) of $C$ and $\ovl{C}$.

\begin{lem}\label{thm1lem}
For $p\in X_d(m,m')$ in $(\ref{fixpath})$, we set $l^{(s)}_i$, $k^{(s)}_i$ and $\delta_i$ as in $(\ref{ldef})$, $(\ref{ldef2})$, $(\ref{kdef})$ and $(\ref{deldef})$. Then we have
\begin{multline}\label{thm1lemclaim}
Q(p)= \prod^{d}_{i=1} \ovl{C}(m-l^{(1)}_i,k^{(1)}_i)\cdots \ovl{C}(m-l^{(\delta_i)}_i,k^{(\delta_i)}_i)\\
\cdot C(m-l^{(\delta_i+1)}_i,|k^{(\delta_i+1)}_i|-1)\cdots C(m-l^{(m-m')}_i,|k^{(m-m')}_i|-1).
\end{multline}
\end{lem}
\nd
{\sl Proof.}

At first, we get $a^{(l^{(\delta_i+1)}_i-1)}_i\leq r$, otherwise, we have $\ovl{r}\leq a^{(l^{(\delta_i+1)}_i-1)}_i$ and hence $l^{(\delta_i)}_i=l^{(\delta_i+1)}_i-1$ and $\ovl{r}\leq a^{(l^{(\delta_i)})}_i=k^{(\delta_i)}_i$ by the definition (\ref{ldef2}) of $l^{(s)}_i$, which contradicts the assumption of $\delta_i$. Further, we get $a^{(l^{(\delta_i+1)}_i-1)}_i=r$ by $a^{(l^{(\delta_i+1)}_i)}_i=k^{(\delta_i+1)}_i\in\{\ovl{j}|1\leq j\leq r\}$ and Definition \ref{pathdef} (iii). Hence we obtain
\begin{equation}\label{thm1lem-rel1}
1\leq a^{(0)}_i\leq a^{(1)}_i\leq\cdots \leq a^{(l^{(\delta_i+1)}_i-1)}_i=r<\ovl{r}\leq a^{(l^{(\delta_i+1)}_i)}_i\leq \cdots\leq a^{(m)}_i\leq\ovl{1}.
\end{equation}

Next, for $0\leq s\leq m-1$ and $1\leq i\leq d$, we set the label $Q(a^{(s)}_i\rightarrow a^{(s+1)}_i)$ as follows:
\begin{equation}\label{thm1lem-pr1}
Q(a^{(s)}_i\rightarrow a^{(s+1)}_i):=
\begin{cases}
\frac{Y_{m-s,a^{(s+1)}_i-1}}{Y_{m-s,a^{(s)}_i}}  & {\rm if}\ 1\leq a^{(s)}_i\leq a^{(s+1)}_i\leq r,\\
\frac{1}{Y_{m-s,|a^{(s+1)}_i|-1}} & {\rm if}\ a^{(s)}_i=r\ {\rm and}\ \ovl{r}\leq a^{(s+1)}_i\leq\ovl{1}, \\
\frac{Y_{m-s,|a^{(s)}_i|}}{Y_{m-s,|a^{(s+1)}_i|-1}} & {\rm if}\ \ovl{r}\leq a^{(s)}_i\leq  a^{(s+1)}_i\leq\ovl{1},
\end{cases}
\end{equation}
which means that the label $Q^{(s)}(p)$ of the edge vt$(m-s;a^{(s)}_1,\cdots,a^{(s)}_d)\rightarrow {\rm vt}(m-s-1;a^{(s+1)}_1,\cdots,a^{(s+1)}_d)$ is as follows (see Definition \ref{labeldef} (i)):
\[ Q^{(s)}(p)=\prod^{d}_{i=1} Q(a^{(s)}_i\rightarrow a^{(s+1)}_i). \]
Therefore, we get
\[ Q(p)=\prod^{m-1}_{s=0}\prod^{d}_{i=1} Q(a^{(s)}_i\rightarrow a^{(s+1)}_i), \]
which is obtained from Definition \ref{labeldef} (ii). To calculate $\prod^{m-1}_{s=0} Q(a^{(s)}_i\rightarrow a^{(s+1)}_i)$ for $1\leq i\leq d$, let us divide the range of product $\prod^{m-1}_{s=0}$ as follows:
\[
\prod^{l^{(\delta_i+1)}_{i}-2}_{s=0},\qq 
\prod^{l^{(m-m')}_i}_{s=l^{(\delta_i+1)}_i-1}\qq {\rm and}\qq
\prod^{m-1}_{s=l^{(m-m')}_i+1},
\]
where in the case $\delta_i=m-m'$, we set
\begin{equation}\label{dmdmdm}
l^{(m-m'+1)}_i:=l^{(m-m')}_i+2.
\end{equation} 
First, let us consider the first range of the product. For $0\leq s\leq l^{(\delta_i+1)}_{i}-2$, using (\ref{thm1lem-rel1}) and (\ref{thm1lem-pr1}), we get
\[ Q(a^{(s)}_i\rightarrow a^{(s+1)}_i)=
\begin{cases}
\frac{Y_{m-s,a^{(s)}_i-1}}{Y_{m-s,a^{(s)}_i}}=\ovl{C}(m-s,a^{(s)}_i)  & {\rm if}\ a^{(s+1)}_i=a^{(s)}_i, \\
1 & {\rm if}\ a^{(s+1)}_i=a^{(s)}_i+1,
\end{cases}
 \]
which means
\begin{equation}\label{thm1lempr2}
\prod^{l^{(\delta_i+1)}_{i}-2}_{s=0}\left(Q(a^{(s)}_i\rightarrow a^{(s+1)}_i)\right)
=\prod^{\delta_i}_{\zeta=1}\ovl{C}(m-l^{(\zeta)}_i,k^{(\zeta)}_i),
\end{equation}
by (\ref{ldef2}) and $k^{(\zeta)}_i:=a^{(l^{(\zeta)}_i)}_i$. 

Next, we consider the second range of the product. If $r-m'\geq i$ then $r\geq m'+i=a^{(m)}_i\geq\cdots\geq a^{(1)}_i\geq a^{(0)}_i$, which implies $\delta_i=m-m'$, and $\prod^{l^{(m-m')}_i}_{s=l^{(\delta_i+1)}_i-1}\left(Q(a^{(s)}_i\rightarrow a^{(s+1)}_i)\right)=1$ by (\ref{dmdmdm}). So we consider the case $r-m'<i$. For $s=l^{(\delta_i+1)}_i-1$, we get
$Q(a^{(s)}_i\rightarrow a^{(s+1)}_i)=\frac{1}{Y_{m-s,|a^{(s+1)}_i|-1}}$, and for $l^{(\delta_i+1)}_i\leq s\leq l^{(m-m')}_i$, $Q(a^{(s)}_i\rightarrow a^{(s+1)}_i)=\frac{Y_{m-s,|a^{(s)}_i|}}{Y_{m-s,|a^{(s+1)}_i|-1}}$ by (\ref{thm1lem-rel1}) and (\ref{thm1lem-pr1}). Thus, we obtain
\begin{eqnarray}\label{thm1lempr3}
& &\prod^{l^{(m-m')}_i}_{s=l^{(\delta_i+1)}_i-1}\left(Q(a^{(s)}_i\rightarrow a^{(s+1)}_i)\right)\nonumber \\
&=&\left(\frac{1}{Y_{m-l^{(\delta_i+1)}_i+1,|a^{(l^{(\delta_i+1)}_i)}_i|-1}}\right)\cdot
\prod^{l^{(m-m')}_i}_{s=l^{(\delta_i+1)}_i} \frac{Y_{m-s,|a^{(s)}_i|}}{Y_{m-s,|a^{(s+1)}_i|-1}}\qq \qq \nonumber\\
&=&\left(\prod^{l^{(m-m')}_i}_{s=l^{(\delta_i+1)}_i} \frac{Y_{m-s,|a^{(s)}_i|}}{Y_{m-s+1,|a^{(s)}_i|-1}}\right)\cdot 
\frac{1}{Y_{m-l^{(m-m')}_i, |a^{(l^{(m-m')}_i+1)}_i|-1}} \nonumber\\
&=&\left(\prod^{m-m'}_{\zeta=\delta_i+1} C(m-l^{(\zeta)}_i,|k^{\zeta}_i|-1)\right)\cdot
\frac{1}{Y_{m'+i-r, d-i}},
\end{eqnarray}
where for the third equality, we used $l^{(m-m')}_i=m-m'-i+r$ (Lemma \ref{kpro2}) and $|a^{(m-m'-i+r+1)}_i|=d-i+1$ (Lemma \ref{cancellem}).

Finally, we consider the last range of the product. Using Lemma \ref{cancellem}, Lemma \ref{kpro2} and (\ref{ldef}), we obtain
\begin{equation}\label{thm1lempr4}
\prod^{m-1}_{s=l^{(m-m')}_i+1}\left(Q(a^{(s)}_i\rightarrow a^{(s+1)}_i)\right)=
\begin{cases}
\prod^{m-1}_{s=m-m'-i+r+1} \frac{Y_{m-s, d-i+1}}{Y_{m-s, d-i}} & {\rm if}\ r-m'<i, \\
1 & {\rm if}\ r-m'\geq i.
\end{cases}
\end{equation}
By (\ref{thm1lempr2}), (\ref{thm1lempr3}) and (\ref{thm1lempr4}), to prove (\ref{thm1lemclaim}), we need to show that
\begin{equation}\label{thm1lempr5}
\prod^d_{i=r-m'+1}\left(\frac{1}{Y_{m'+i-r, d-i}} \prod^{m-1}_{s=m-m'-i+r+1}\frac{Y_{m-s, d-i+1}}{Y_{m-s, d-i}}\right)=1. 
\end{equation}
We set 
\[A:=\prod^d_{i=r-m'+1}\left(\frac{1}{Y_{m'+i-r, d-i}}\right),\ {\rm and}\q B:=\prod^d_{i=r-m'+1}\left(\prod^{m-1}_{s=m-m'-i+r+1}\frac{Y_{m-s, d-i+1}}{Y_{m-s, d-i}}\right).\]

We obtain the followings:
\[
A=\prod^d_{i=r-m'+1}\left(\frac{1}{Y_{m'+i-r, d-i}}\right)=
\prod^{d-1}_{i=r-m'+1}\left(\frac{1}{Y_{m'+i-r, d-i}}\right)=
\prod^{m'+d-r-1}_{k=1}\left(\frac{1}{Y_{k, d-r+m'-k}}\right),
\]
and
\begin{eqnarray*}
B&=&\prod^d_{i=r-m'+1}\left(\prod^{m-1}_{s=m-m'-i+r+1}\frac{Y_{m-s, d-i+1}}{Y_{m-s, d-i}}\right)=
\prod^d_{i=r-m'+1}\left(\prod^{m'+i-r-1}_{s=1}\frac{Y_{s, d-i+1}}{Y_{s, d-i}}\right)\\
&=&\prod^{m'+d-r-1}_{s=1}\left(\frac{Y_{s,d-r+m'-s}}{Y_{s,d-r+m'-s-1}}\frac{Y_{s,d-r+m'-s-1}}{Y_{s,d-r+m'-s-2}}
\frac{Y_{s,d-r+m'-s-2}}{Y_{s,d-r+m'-s-3}}\cdots \frac{Y_{s,1}}{Y_{s,0}} \right)\\
&=&\prod^{m'+d-r-1}_{s=1} Y_{s,d-r+m'-s},
\end{eqnarray*}
where  note that $Y_{s,0}=1$ (see Remark \ref{importantrem}). Thus we have $A\cdot B=1$,
which implies $(\ref{thm1lempr5})$. \qed

\vspace{3mm}

Let us prove the main theorem.

\nd
{\sl Proof of Theorem \ref{thm1}.}

Using Lemma \ref{thm1lem}, we see that $Q(p)$ $(p\in X_d(m,m'))$ is described as (\ref{thm1lemclaim}) with $\{k^{(s)}_i\}_{1\leq i\leq d, 1\leq s\leq m-m'}$ which satisfy the conditions in Lemma \ref{kpro2} (ii), that is, $1\leq k^{(s)}_1<k^{(s)}_2<\cdots<k^{(s)}_d\leq\ovl{1}$. If $m'+i\leq r$, then $a^{(0)}_i\leq a^{(1)}_i\leq \cdots\leq a^{(m)}_i=m'+i$, which means that $1\leq k^{(1)}_i\leq k^{(2)}_i\leq\cdots\leq k^{(m-m')}_i\leq m'+i$ for $1\leq i\leq r-m'$. For $r-m'+1\leq i\leq d$, the inequality (\ref{iseqine}) implies $1\leq k^{(1)}_i\leq k^{(2)}_i\leq\cdots\leq k^{(m-m')}_i\leq \ovl{1}$. Thus, $\{k^{(s)}_i\}$ satisfies the conditions $(*)$ in Theorem \ref{thm1}.

Conversely, let $\{K^{(s)}_i\}_{1\leq i\leq d, 1\leq s\leq m-m'}$ the set of numbers which satisfies the conditions $(*)$ in Theorem \ref{thm1}:
\begin{equation}\label{cond1}
1\leq K^{(s)}_1<K^{(s)}_2<\cdots<K^{(s)}_d\leq \ovl{1}\q (1\leq s\leq m-m'),
\end{equation}
\begin{equation}\label{cond2}
1\leq K^{(1)}_i\leq \cdots\leq K^{(m-m')}_i\leq m'+i\q (1\leq i\leq r-m'),
\end{equation}
and
\begin{equation}\label{cond3}
1\leq K^{(1)}_i\leq \cdots\leq K^{(m-m')}_i\leq \ovl{1}\q (r-m'+1\leq i\leq d).
\end{equation}
We need to show that there exists a path $p\in X_d(m,m')$ such that
\begin{multline}\label{finclaim}
Q(p)= \prod^{d}_{i=1} \ovl{C}(m-L^{(1)}_i,K^{(1)}_i)\cdots \ovl{C}(m-L^{(\delta_i)}_i,K^{(\delta_i)}_i)\\
\cdot C(m-L^{(\delta_i+1)}_i,|K^{(\delta_i+1)}_i|-1)\cdots C(m-L^{(m-m')}_i,|K^{(m-m')}_i|-1),
\end{multline}
where $\delta_i$ $(1\leq\delta_i\leq m-m')$ are the numbers which satisfy $1\leq K^{(1)}_i\leq \cdots\leq K^{(\delta_i)}_i\leq r<\ovl{r}\leq K^{(\delta_i+1)}_i\leq \cdots \leq K^{(m-m')}_i\leq\ovl{1}$, and
\[ L^{(s)}_{i}:=
\begin{cases}
K^{(s)}_{i}+s-i-1 & {\rm if}\  K^{(s)}_{i}\in\{j|1\leq j\leq r\}, \\
s-i+r & {\rm if}\ K^{(s)}_{i}\in\{\ovl{j}|1\leq j\leq r\},
\end{cases}
\]
for $1\leq s\leq m-m'$ and $1\leq i\leq d$. Since we supposed $K^{(s)}_i<K^{(s)}_{i+1}$, we can easily verify \begin{equation}\label{lequ1}
L^{(s)}_i\leq L^{(s)}_{i+1}\q {\rm if}\ K^{(s)}_i\in\{j|1\leq j\leq r\},
\end{equation}
and 
\begin{equation}\label{lequ2}
L^{(s)}_i=L^{(s)}_{i+1}+1\q {\rm if}\ K^{(s)}_i\in\{\ovl{j}|1\leq j\leq r\}. 
\end{equation}
We claim that $0\leq L^{(s)}_i\leq m-1$. By the condition (\ref{cond1}), we get $i\leq K^{(s)}_i$. So it is clear that $0\leq L^{(s)}_i$. For $1\leq i\leq r-m'$ and $1\leq s\leq m-m'$, it follows from the condition $(\ref{cond2})$ that $L^{(s)}_{i}=K^{(s)}_{i}+s-i-1\leq m'+i+s-i-1=m'+s-1\leq m-1$. For $r-m'<i$, we get $L^{(s)}_i\leq r-i+s<m'+s\leq m$. Therefore, we have $0\leq L^{(s)}_i\leq m-1$ for all $1\leq i\leq d$ and $1\leq s\leq m-m'$.

Note that if $K^{(s)}_i\in\{\ovl{j}|1\leq j\leq r\}$, then $\ovl{r}\leq K^{(s)}_i<K^{(s)}_{i+1}\leq \ovl{1}$ and hence
\begin{equation}\label{lequ3}
L^{(s+1)}_i=L^{(s)}_i+1.
\end{equation}
We define a path $p={\rm vt}(m;a^{(0)}_1,\cdots,a^{(0)}_d)\rightarrow\cdots \rightarrow {\rm vt}(0;a^{(m)}_1,\cdots,a^{(m)}_d)\in X_d(m,m')$ as follows: For $i$ $(1\leq i\leq r-m')$, we define the $i$-sequence (Definition \ref{iseq}) of $p$ as

\begin{flushleft}
$a^{(0)}_i=i,\ a^{(1)}_i=i+1,\ a^{(2)}_i=i+2,\cdots,a^{(L^{(1)}_i)}_i=i+L^{(1)}_i$,
 
$ a^{(L^{(1)}_i+1)}_i=i+L^{(1)}_i,\ a^{(L^{(1)}_i+2)}_i=i+L^{(1)}_i+1,\cdots, a^{(L^{(2)}_i)}_i=i+L^{(2)}_i-1, $

$ a^{(L^{(2)}_i+1)}_i=i+L^{(2)}_i-1,\ a^{(L^{(2)}_i+2)}_i=i+L^{(2)}_i,\cdots, a^{(L^{(3)}_i)}_i=i+L^{(3)}_i-2, $
\end{flushleft}
\begin{equation}\label{jlist2} \vdots \end{equation}
\begin{flushleft}
$a^{(L^{(m-m'-1)}_i+1)}_i=i+L^{(m-m'-1)}_i-m+m'+2,\ 
\cdots, a^{(L^{(m-m')}_i)}_i=i+L^{(m-m')}_i-m+m'+1,$

$a^{(L^{(m-m')}_i+1)}_i=i+L^{(m-m')}_i-m+m'+1,\ a^{(L^{(m-m')}_i+2)}_i=i+L^{(m-m')}_i-m+m'+2,$

$a^{(L^{(m-m')}_i+3)}_i=i+L^{(m-m')}_i-m+m'+3,\cdots, a^{(m)}_i=m'+i$.
\end{flushleft}

\vspace{2mm}

For $i$ $(r-m'+1\leq i\leq d)$, we define the $i$-sequence of $p$ as

\begin{flushleft}
$a^{(0)}_i=i,\ a^{(1)}_i=i+1,\ a^{(2)}_i=i+2,\cdots,a^{(L^{(1)}_i)}_i=i+L^{(1)}_i$,
 
$ a^{(L^{(1)}_i+1)}_i=i+L^{(1)}_i,\ a^{(L^{(1)}_i+2)}_i=i+L^{(1)}_i+1,\cdots, a^{(L^{(2)}_i)}_i=i+L^{(2)}_i-1, $

$ a^{(L^{(2)}_i+1)}_i=i+L^{(2)}_i-1,\ a^{(L^{(2)}_i+2)}_i=i+L^{(2)}_i,\cdots, a^{(L^{(3)}_i)}_i=i+L^{(3)}_i-2, $
\end{flushleft}
\begin{equation}\label{jlist3} \vdots \end{equation}
\begin{flushleft}
$a^{(L^{(\delta_i-1)}_i+1)}_i=i+L^{(\delta_i-1)}_i-\delta_i+2,\ 
\cdots, a^{(L^{(\delta_i)}_i)}_i=i+L^{(\delta_i)}_i-\delta_i+1,$

$a^{(L^{(\delta_i)}_i+1)}_i=i+L^{(\delta_i)}_i-\delta_i+1,\ a^{(L^{(\delta_i)}_i+2)}_i=i+L^{(\delta_i)}_i-\delta_i+2,$

$a^{(L^{(\delta_i)}_i+3)}_i=i+L^{(\delta_i)}_i-\delta_i+3,\cdots,\ a^{(L^{(\delta_i+1)}_i-1)}_i=r$,

$a^{(L^{(\delta_i+1)}_i)}_i=K^{(\delta_i+1)}_i,\ a^{(L^{(\delta_i+2)}_i)}_i=K^{(\delta_i+2)}_i,\cdots,\ a^{(L^{(m-m')}_i)}_i=K^{(m-m')}_i$, 

$a^{(L^{(m-m')}_i+1)}_i=a^{(L^{(m-m')}_i+2)}_i=\cdots=a^{(m)}_i=\ovl{d-i+1}$.
\end{flushleft}

\vspace{2mm}

It is easy to see that $a^{(L^{(s)}_i)}_i=K^{(s)}_i$ $(1\leq s\leq m-m')$ by the above lists. Clearly, the path $p$ satisfies Definition \ref{pathdef} (iii) and (iv). For $1\leq s\leq L^{(\delta_{i}+1)}_i-1$, we obtain $a^{(s)}_i<a^{(s)}_{i+1}$ by (\ref{lequ1}). For $\delta_{i}+1\leq s\leq m-m'$, we obtain $a^{(L^{(s)}_i)}_i<a^{(L^{(s)}_{i})}_{i+1}$ since  $a^{(L^{(s)}_i)}_i=K^{(s)}_i<K^{(s)}_{i+1}=a^{(L^{(s)}_{i+1})}_{i+1}=a^{(L^{(s)}_{i}-1)}_{i+1}\leq
a^{(L^{(s)}_{i})}_{i+1}$ by (\ref{cond1}) and (\ref{lequ2}). For $L^{(m-m')}_i+1\leq s\leq m$, we obtain $a^{(s)}_i=\ovl{d-i+1}$, and then we get $a^{(s)}_{i+1}=\ovl{d-i}$ since $L^{(m-m')}_{i+1}=L^{(m-m')}_i-1<L^{(m-m')}_i\leq s$, which means $a^{(s)}_i<a^{(s)}_{i+1}$. Therefore, $a^{(s)}_i<a^{(s)}_{i+1}$ for all $1\leq i\leq d-1$ and $1\leq s\leq m-m'$, which means the path $p$ satisfies Definition \ref{pathdef} (ii).

Finally, for $a^{(s)}_i\in \{\ovl{j}|1\leq j\leq r\}$, we need to verify $a^{(s)}_i<a^{(s-1)}_{i+1}$. The definition (\ref{jlist2}), (\ref{jlist3}) of $i$-sequence of $p$ shows that either $s=L^{(\zeta)}_i$ for some $\zeta$ $(\delta_i+1\leq\zeta\leq m-m')$ or $L^{(m-m')}_i<s$. In the case $s=L^{(\zeta)}_i$, using (\ref{cond1}) and (\ref{lequ2}), we see that $a^{(s)}_i=a^{(L^{(\zeta)}_i)}_i=K^{(\zeta)}_i<K^{(\zeta)}_{i+1}=
a^{(L^{(\zeta)}_{i+1})}_{i+1}=a^{(L^{(\zeta)}_{i}-1)}_{i+1}=a^{(s-1)}_{i+1}$. In the case $L^{(m-m')}_i<s$, we obtain $a^{(s)}_i=\ovl{d-i+1}<\ovl{d-i}=a^{(s-1)}_{i+1}$ since $L^{(m-m')}_{i+1}=L^{(m-m')}_{i}-1<s-1$. Therefore, we have $a^{(s)}_i<a^{(s-1)}_{i+1}$ for $a^{(s)}_i\in \{\ovl{j}|1\leq j\leq r\}$, which means the path $p$ satisfies Definition \ref{pathdef} (v).

Hence $p$ is well-defined, and (\ref{finclaim}) is follows from Lemma \ref{thm1lem}, and Theorem \ref{thm1} follows from Proposition \ref{pathlem}. \qed


\end{document}